\documentclass[xypic,syntonly,amssymb,verbatim,12pt]{amsart}
\usepackage{graphicx}
\usepackage[all]{xy}
\usepackage{amssymb}
\usepackage{amsmath}
\usepackage[mathscr]{euscript}

\usepackage{epsfig}
\theoremstyle{plain}
\newtheorem{Thm}{Theorem}
\newtheorem{Prop}{Proposition}
\newtheorem{Cor}{Corollary}
\newtheorem{Lem}{Lemma}
\newtheorem{Def}{Definition}
 \theoremstyle{definition}
\theoremstyle{remark}
\newtheorem{Rems} {Remark}

\numberwithin{equation}{section}
\providecommand{\norm}[1]{\lVert#1\rVert}

\begin{document}
  \title {Holomorphic Yang-Mills fields on $B$-branes}

 \author{ Andr\'{e}s   Vi\~{n}a} 
\address{Departamento de F\'{i}sica. Universidad de Oviedo.    Garc\'{\i}a Lorca 18.
     33007 Oviedo. Spain. }
\email{vina@uniovi.es}

 \maketitle
\begin{abstract}
 Considering $B$-branes over a complex manifold $X$ as objects of the bounded derived category of coherent sheaves over $X$, we define holomorphic gauge fields on $B$-branes and introduce the Yang-Mills functional for these fields. These definitions extend well-known concepts in the context of vector bundles to the setting of $B$-branes.

For a given $B$-brane, we show that its Atiyah class is the obstruction to
the existence of gauge fields. When $X$ is the variety of complete flags in a $3$-dimensional complex vector space, we prove that any $B$-brane over $X$ admits at most one holomorphic gauge field.

Furthermore, we establish that the set of Yang-Mills fields on a given $B$-brane, if nonempty, is in bijective correspondence with the points of an algebraic set defined by $m$ complex polynomials of degree less than four in $m$ indeterminates, where $m$ is the dimension of the space of morphisms from the brane to its tensor product with the sheaf of holomorphic one-forms.

\end{abstract}

\section {Introduction} \label{S:intro}

  In this article, we extend  the well-known concepts of gauge field and Yang-Mills field on vector bundles to $B$-branes.
From the mathematical point of view, a $B$-brane over  a complex manifold $X$ is an object of $D^b(X)$, the bounded derived category of  coherent analytic sheaves over $X$
 \cite[Sect. 5.4]{Aspin} \cite[Sect. 5.3]{Aspin-et}. 

The simplest $B$-branes are the  {\em holomorphic } vector bundles.  A gauge field on a holomorphic vector bundle $V\to M,$ in mathematical terms a connection on $V,$  defines  a derivation of 
sections of $V$ along vector fields on $M$  giving rise to   sections 
not necessarily holomorphic. A holomorphic gauge field  on $V,$ i. e. a holomorphic connection on $V$ \cite{Atiyah}, allows us
 to define a derivative of the holomorphic sections of $V$ along any ``direction'', giving rise to {\em holomorphic} sections.  
 
 Not every holomorphic vector bundle admits a holomorphic connection. The obstruction for this existence is the Atiyah class
of the bundle. 
The vanishing of this class is a necessary and sufficient condition for the existence of such connections
 on the 
vector bundle.
 Assuming that there is no obstruction to the existence of holomorphic connections, the set of these connections on a 
 vector bundle over a compact manifold is a {\it finite dimensional} affine space, in contrast to the smooth case. This property is a consequence of the finite dimensionality of the coherent cohomology groups.
 
The two features discussed above concerning the holomorphic connections on vector bundles: the   obstruction to its existence,  and the finite dimension of the   space of those fields, when it is not empty, are also present in the extension 
to general branes.


{\it Connections on coherent sheaves.}	
Given a coherent sheaf ${\mathscr F}$   over the compact analytic manifold $X,$
a holomorphic connection on ${\mathscr F}$ \cite{Katz} determines 
   isomorphisms between the stalks of ${\mathscr F}$ over ``close'' points of $X;$   i. e.   identifications in the directions defined by the tangent vectors of the base.
 
  The idea of being
	``infinitesimally close'' can be formulated by means of  the first infinitesimal neighborhood of the diagonal of $X$. In this setting, a  connection  
	on ${\mathscr F}$ can be regarded 
	as a right inverse of the natural morphism ${\mathscr J}^1({\mathscr F})\to{\mathscr F}$, where ${\mathscr J}^1({\mathscr F})$  is the corresponding  $1$-jet sheaf \cite[Sect.\,3]{Vina21}. This inverse determines a morphism of abelian sheaves  
	$\nabla:{\mathscr F}\to \Omega^1({\mathscr F})=\Omega^1_X\otimes_{{\mathscr O}_X}{\mathscr F},$ which satisfies the Leibniz rule.

	In this context, the obstruction to the existence of a holomorphic connection on the sheaf ${\mathscr F}$ is 
	 an element 
	  of the group ${\rm Ext}^1({\mathscr F},\,\Omega^1({\mathscr F})).$ Furthermore, when the set of holomorphic connections 
	 on ${\mathscr F}$ is nonempty, it is an  affine space   associated  to the {\it finite} dimensional vector space  
	 ${\rm Hom}_{{\mathscr O}_X}({\mathscr F},\,\Omega^1({\mathscr F}))$.


{\it Yang-Mills fields on sheaves.} The reflexive sheaves 
	 might be thought as ``vector bundles with singularities'' 
	 \cite[p.\,121]{Hartshorne1}, and these singularities    can in some way be avoided.
	Setting ${\mathcal S}$ for the singularity set of the reflexive sheaf ${\mathscr F}$, let us assume that on the locally free sheaf ${\mathscr F}|_{X\setminus {\mathcal S}}$ there is defined a Hermitian metric, we say that ${\mathscr F}$ is a Hermitian sheaf. If $X$ is a K\"ahler manifold, on  the set of holomorphic connections over ${\mathscr F},$ one defines the Yang-Mills functional $\mathcal{YM}$. Essentially, the value  of  $\mathcal{YM}$ at a connection $\nabla$
	is the squared norm 
	$\norm{K_{\nabla}}^2$ of the curvature of that connection 
	\cite[p.\,417]{Hamilton}, \cite[p.\,44]{Moore}, \cite[p.\,357]{Naber}. 
	
	  The stationary points of  $\mathcal{YM}$ are the Yang-Mills fields on ${\mathscr F}$. The set of these   points of this functional will be denoted by ${\sf YM}({\mathscr F})$.
	 Denoting 
	 $m:={\rm dim}\,{\rm Hom}_{{\mathscr O}_X}({\mathscr F},\,\Omega^1({\mathscr F})),$
  we will prove 
that if the reflexive sheaf ${\mathscr F}$ admits holomorphic
connections,
 then the set ${\sf YM}({\mathscr F})$  is in bijective correspondence with the points of an algebraic subset of ${\mathbb C}^m$ defined by $m$ algebraic equations of degree $\leq 3$ (Theorem \ref{P:numeroYM}). Thus, in general, the set of Yang-Mills fields on ${\mathscr F}$ will be a finite set.


{\it Yang-Mills fields on a $B$-brane.}	
	As we said, a $B$-brane on $X$ is 
	a complex   $({\mathscr F}^{\bullet},\,\delta^{\bullet})$ of analytic coherent sheaves on $X.$  According to the preceding paragraph, it is reasonable to define a holomorphic gauge field on this brane as an element  
 $\psi\in {\rm Hom}_{D^b(X)}\big({\mathscr F}^{\bullet},\,  {\mathscr J}^1({\mathscr F}^{\bullet})\big)$ which lifts the identity on  ${\mathscr F}^{\bullet}.$

When the manifold $X$ admits a stratification, where the strata satisfy certain properties, the space of holomorphic gauge fields on any $B$-brane on $X$ is a set with cardinal $\leq 1$ (Theorem \ref{Th:Z_i}). The case where $X$ is the projective space ${\mathbb P}^n$ has been considered in \cite{Vina24}. Here, we study the case where $X$ is the variety of complete flags in ${\mathbb C}^3$ (Theorem \ref{Thm:flag}).
 
 The homomorphic gauge field $\psi$ determines a unique morphism between each of the cohomology sheaves
  $\vartheta^j:{\mathscr H}^j({\mathscr F}^{\bullet})\to {\mathscr H}^j(\Omega^1({\mathscr F}^{\bullet})),$ which in fact is a holomorphic connection   ${\mathscr H}^j({\mathscr F}^{\bullet})$.
  
   When the cohomology sheaves ${\mathscr H}^j$ are Hermitian, we define the value of the Yang-Mills functional on the above  gauge field $\psi$ as
 $\sum_i(-1)^i\norm{K_{\vartheta_i}}^2.$ Thus, the Yang-Mills functional is a kind of Euler characteristic  of the gauge field. 
 Obviously, this definition of the Yang-Mills functional 
 in the context of $B$-branes   generalizes the one  for coherent sheaves.

The gauge field $\psi$ is a Yang-Mills field if it is a stationary point of Yang-Mills functional. 
In Theorem \ref{P:numero YM}, we generalize the result given in Theorem \ref{P:numeroYM}, about the cardinal of the set of Yang-Mills fields on sheaves, to a general brane ${\mathscr F}^{\bullet}$.
 We will prove that if the $\vartheta^j$ are Yang-Mills fields on the ${\mathscr H}^j$,  so is $\psi$  on the brane
  (Proposition \ref{estationaryY-Mfunctional}). Theorem \ref{Th:inverse} is a partial converse to Proposition \ref{estationaryY-Mfunctional}. 
 
 

The article is organized in two sections. In Section \ref{S:Holomorphi}
 are considered the holomorphic gauge fields on  $B$-branes. Subsection \ref{ss:FirstJet} concerns with the holomorphic connections on a sheaf. In Subsection \ref{Ss:Fields_brane}, we define the holomorphic gauge fields on a brane and discuss the obstruction to their existence. 
  Subsection \ref{Ss:Flag} discusses branes on a manifold which admits certain stratifications; in this subsection we prove Theorems \ref{Th:Z_i} and \ref{Thm:flag}.
  
   Section \ref{ss:holomorphicYang} is devoted to the holomophic Yang-Mills fields. In Subsection \ref{Ss:Hermitian}, we introduce the Yang-Mills functional for holomorphic connections on reflexive sheaves and prove Theorem  \ref{P:numeroYM}.
  The Yang-Mills functional for holomorphic gauge fields on a $B$-brane is defined in Subsection \ref{Ss:YM-Brenes}, along with a justification for this definition. 
   In this subsection are also proved Proposition \ref{estationaryY-Mfunctional} and Theorems   \ref{P:numero YM} and \ref{Th:inverse}.

\section{Holomorphic gauge fields on $B$-branes}\label{S:Holomorphi}
As we mentioned in the Introduction, our purpose is to define holomorphic
gauge fields on $B$-branes, extending the concept of holomorphic connection
on vector bundles.


\subsection{Holomorphic connections on a sheaf.} \label{ss:FirstJet} 
 The definition of connection on a coherent sheaf, as it is introduced 
   for example in \cite{Katz}, is reformulated in this section, so that it is suitable for translation to objects in the derived category of coherent sheaves.

As explained in the Introduction, given a coherent sheaf ${\mathscr F}$ over the complex manifold $X$,
 the existence of a  homomorphic connection on  
 ${\mathscr F}$ 
  should define an isomorphism between the stalks  of 
 the sheaf at any two ``infinitesimally close'' points of $X.$
By means of the first infinitesimal  neighborhood $X^{(1)}$ of  diagonal of $X$ \cite[p.\,698]{G-H}, it is possible to formulate this idea in a suitable way also for algebraic varieties.

  The outline of a new definition for a connection on a sheaf ${\mathscr G}$ over a complex variety $X$, inspired by an idea from crystalline cohomology, is as follows.
  If $R$ is a ${\mathbb C}$-algebra, ${\rm Hom}({\rm Spec}\, R,\,X)$ is the set of points of complex algebraic variety $X$ with values  in $R$. Two points $x,y$ are infinitesimally close if the morphism $(x,\,y):{\rm Spec}\, R\to X\times X$ factorizes through the infinitesimal neighborhood $X^{(1)}.$  
 Hence, a connection on the sheaf ${\mathscr G }$ should define an identification of the pullbacks $x^*{\mathscr G}$ and
  $y^*{\mathscr G}$ for any two infinitesimally close points.
  
   Therefore, following Deligne \cite[p.\,6]{Deligne}, one defines a connection   on the coherent ${\mathscr O}_X$-module ${\mathscr F}$ as an element of 
\begin{equation}\label{Deligne_Gauge}
{\rm Hom}_{X^{(1)}}(\pi_1^*{\mathscr F},\, \pi_2^*{\mathscr F}),
\end{equation}
 which is the identity on $X$, where $\pi_1,\pi_2:X^{(1)}\to X$ are the projections and $\pi_i^*$ the corresponding inverse image functor.
 By the adjunction isomorphism, a  holomorphic connection  on ${\mathscr F}$ is a morphism ${\mathscr F} \to\pi_{1*} \pi^*_2{\mathscr F}$ which is a right inverse of the projection
 $\pi:\pi_{1*} \pi^*_2{\mathscr F}\to{\mathscr F}.$ 

As  $\pi_{1*} \pi^*_2{\mathscr F}$ is
 the first jet sheaf   of the coherent sheaf ${\mathscr F}$, 
in the following paragraph we review the definition of the  jet sheaf. 

 \subsubsection{The first jet sheaf}\label{Ss: first jet}

Let  $X$ be an analytic complex compact connected manifold. 
 Let   $i:\Delta\hookrightarrow X\times X$ be the embedding of the diagonal. As a closed submanifold, $\Delta$ is defined by an ideal $\mathscr I$ of
 $\Hat{\mathscr O}:={\mathscr O}_{X\times X}.$
 The first infinitesimal neighborhood of $\Delta$ is the following ringed space
  $$X^{(1)}=\Big(\Delta,\,{\mathscr O}_{X^{(1)}}:=\big(\Hat{\mathscr O}/{\mathscr I}^2\big)|_{\Delta}\Big)$$

We set $p_1,p_2: X\times X\rightrightarrows X$ for the corresponding projection morphisms. For $a=1,2$, the compositions $p_a\circ i,$
 $ \xymatrix{
 \Delta \ar[r]^i &
 X\times X\ar@/^/[r]^{p_1}\ar@/_/[r]_{p_2}    & X
 } $ will be denoted by $\pi_a.$

 Given ${\mathscr F}$ a left ${\mathscr O}_X$-module,   its inverse image by $\pi_2$ is the left ${\mathscr O}_{X^{(1)}}$-module
  $$\pi_2^*({\mathscr F})={\mathscr O}_{X^{(1)}}\otimes_{\pi_2^{-1}{{\mathscr O}_X}}\pi_2^{-1}{\mathscr F}.$$
And the first jet sheaf ${\mathscr J}^1({\mathscr F})$ of ${\mathscr F}$ is the left ${\mathscr O}$-module defined by (see \cite[p.\,505]{Malgrange} \cite[Sect.\,2.4]{Vakil})
\begin{equation}\label{jet1bundle}
{\mathscr J}^1({\mathscr F})={\pi_1}_*\pi_2^*({\mathscr F}).
\end{equation}

We set $\Omega^1({\mathscr F}):=\Omega^1_X\otimes_{{\mathscr O}_X}{\mathscr F},$ where $\Omega^1_X$ is the sheaf of holomorphic $1$-forms on $X$. The first jet sheaf  ${\mathscr J}^1({\mathscr F})$ is the {\em abelian} sheaf ${\mathscr F}\oplus \Omega^1({\mathscr F})$ endowed with the following left ${\mathscr O}_X$-module structure
\begin{equation}\label{productbyf}
 f\cdot(\sigma\oplus\beta)=f\sigma\oplus(f\beta+df\otimes\sigma).
\end{equation}
One has the Atiyah exact sequence of ${\mathscr O}_X$-modules
\begin{equation}\label{AtiyahExact}
0\to \Omega^1({\mathscr F})\to {\mathscr J}^1({\mathscr F})\overset{\pi}{\rightarrow}{\mathscr F}\to 0,
\end{equation} 
where $\pi$ is the projection morphism.

Since $\pi_1^*$ is the left adjoint of ${\pi_1}_*$, one has
\begin{equation}\label{adjuntion}
{\rm Hom}\,_{{\mathscr O}_X}\big({\mathscr F},\,{\mathscr J}^1({\mathscr F})\big)= {\rm Hom}\,_{{\mathscr O}_{X^{(1)}}} \big(\pi_1^*{\mathscr F}, \,\pi_2^*{\mathscr F}\big).
\end{equation}

We  introduce  a piece of notation  that will be used later on. Let ${\sf Forg}$ denote the forgetful functor from the category of ${\mathscr O}_X$-modules to the category of abelian sheaves over $X.$  We set 
  $$\pmb{\mathscr J}^1({\mathscr F}):={\sf Forg}({\mathscr F)}\oplus{\sf Forg}(\Omega^1({\mathscr F}) ).$$
 In this way ${\mathscr J}^1({\mathscr F})$ is the abelian sheaf $\pmb{\mathscr J}^1({\mathscr F})$ endowed with the ${\mathscr O}_X$-module structure defined in  (\ref{productbyf}). In short, we will write
 $\,{\mathscr J}^1({\mathscr F})={\mathscr F}\,\tilde\oplus\,\Omega^1({\mathscr F}). $

One has the morphism of ${\mathbb C}_X$-modules 
 \begin{equation}\label{eta}
 \eta:{\mathscr F}\to \pmb{\mathscr J}^1({\mathscr F})= {\mathscr F}\oplus\Omega^1({\mathscr F}), 
\;\; \sigma\mapsto \sigma\oplus 0.
\end{equation}
And from (\ref{productbyf}), it follows 
\begin{equation}\label{etadef}
f\eta(\sigma)= \eta(f\sigma)+df\otimes\sigma,
\end{equation}
where $f\in\mathscr{O}_X$ is a function of the structure sheaf of $X$
and $\sigma$ a holomorphic section of ${\mathscr F}.$

\smallskip

 
Given $\phi\in{\rm Hom}_{\mathscr{O}_X}({\mathscr F},\,{\mathscr J}^1({\mathscr F}))$ a right inverse of $\pi$;
that is, such that $\pi\circ\phi={\rm id}$. Then $p\circ(\phi-\eta)=0,$ where $p:\pmb{\mathscr J}^1({\mathscr F})\to {\mathscr F}$ is the projection morphism.
Thus $\phi-\eta$ factors uniquely through ${\rm Ker}(p)=\Omega^1({\mathscr F})$,
defining the morphism $\nabla$ in the following diagram in the category of ${\mathbb C}_X$-modules.
\[
\xymatrix{0 \ar[r] & \Omega^1({\mathscr F})\ar[r] & \pmb{\mathscr J}^1({\mathscr F})\ar[r]^p& {\mathscr F} \ar[r] & 0 \\
 & & & {\mathscr F} \ar[ul]_{\phi-\eta} \ar@{-->}[ull]^{\nabla} & }
\]

From (\ref{etadef}), together with the fact that $\phi$ is an ${\mathscr O}_X$-morphism, it follows that 
  $\nabla$ satisfies   the Leibniz rule $\nabla(f\sigma)=\partial f\otimes\sigma+f\nabla(\sigma).$
Therefore, one can give a new definition of holomorphic connection equivalent to the one given in \cite{Katz}.

\begin{Def}\label{Def2}
The holomorphic connections on the coherent ${\mathscr O}_X$-module ${\mathscr F}$   are the elements of the following set
\begin{equation}\label{Def:connection}
\{\phi\in{\rm Hom}_{\mathscr{O}_X}({\mathscr F},\,{\mathscr J}^1({\mathscr F}))\;|\;  \pi\circ\phi={\rm id}\}.
\end{equation}
\end{Def}

That is,  a holomorphic connection $\phi$ is a splitting of the corresponding Atiyah exact sequence. On the other hand, by (\ref{adjuntion}), $\phi$
 is an element of ${\rm Hom}_{{\mathscr O}_{X^{(1)}}}\big(\pi_1^*{\mathscr F},\,\pi_2^*{\mathscr F}\big)$ which is the identity on $X.$ Thus, we recover   (\ref{Deligne_Gauge}).

\newcommand\Ext{\mathrm{Ext}}

\subsection{Gauge fields on a $B$-brane} \label{Ss:Fields_brane}  Let $({\mathscr F}^{\bullet},\,\delta^{\bullet})$ be a $B$-brane  on the complex manifold $X$; that is, ${\mathscr F}^{\bullet}$ is an object of the   category $D^b(X)$, the bounded derived category of coherent sheaves over $X$. The corresponding first jet complex ${\mathscr J}^1({\mathscr F}^{\bullet})$ is defined by
$${\mathscr J}^1({\mathscr F}^{\bullet}):=R\pi_{1*} L\pi_2^*{\mathscr F}^{\bullet}\simeq{\mathscr O}_{X^{(1)}}\otimes^L{\mathscr F}^{\bullet}.$$ 

As ${\mathscr O}_{X^{(1)}}$ is the locally free module ${\mathscr O}_X\oplus \Omega^1_X$, then ${\mathscr J}^1({\mathscr F}^{\bullet})$ is the complex of abelian sheaves
 $\sf{Forg}({\mathscr F}^{\bullet})\,\oplus \,\sf{Forg}(\Omega^1({\mathscr F}^{\bullet})),$
  with ${\mathscr O}_X$-structure given by
\begin{equation}\label{fsigmabullet}
f\cdot(\sigma^{\bullet}\oplus \beta^{\bullet})=f\sigma^{\bullet}\oplus(df\otimes\sigma^{\bullet}+f\beta^{\bullet}).
\end{equation}

More precisely, one can consider the complex of abelian sheaves direct sum of ${\mathscr F}^{\bullet}$
and $\Omega^1({\mathscr F}^{\bullet}).$ According to the notation introduced at the end of Subsection \ref{Ss: first jet}, we write
$$\,\pmb{\mathscr J}^1({\mathscr F}^{\bullet}):= \big({\sf Forg}(
{\mathscr F}^{\bullet})\oplus {\sf Forg}( \Omega^1({\mathscr F}^{\bullet})),\,\,\delta_J^{\bullet}\,\big),$$
where $\delta_J^{\bullet}=\delta^{\bullet}\oplus(1\otimes\delta^{\bullet}).$
The complex $\pmb{\mathscr J}^1({\mathscr F}^{\bullet})$ can be equipped with an ${\mathscr O}_X$-module structure  (\ref{fsigmabullet}).
As   $\delta$ is the coboundary   operator of a complex ${\mathscr O}_X$-modules, then $\delta_J(f\cdot(-))=f\cdot(\delta_J(-)).$ Thus, ${\mathscr J}^1({\mathscr F}^{\bullet})$ is the complex $\pmb{\mathscr J}^1({\mathscr F}^{\bullet})$  endowed with this ${\mathscr O}_X$-module structure (\ref{fsigmabullet}).
 
One  has the exact sequence of complexes of ${\mathscr O}_X$-modules 
\begin{equation}\label{ExactJ1}
0\to\Omega^1({\mathscr F}^{\bullet})\overset{\iota}{\rightarrow}{\mathscr J}^1({\mathscr F}^{\bullet} ) \overset{\pi}{\rightarrow} {\mathscr F}^{\bullet}\to 0.
\end{equation}
The morphisms in the derived category $D^b(X)$ defined by the arrows in the sequence (\ref{ExactJ1}) will  also be denoted by $\iota$ and $\pi.$

%

According  to Definition \ref{Def2}, we propose the following definition.

\begin{Def}\label{Def-brane} A gauge field on ${\mathscr F}^{\bullet}$ is an element $\psi\in{\rm Hom}_{D^b(X)}({\mathscr F}^{\bullet},\, {\mathscr J}^1({\mathscr F}^{\bullet}))$, such that $ {\pi\circ}\psi ={\rm id}\in {\rm Hom}_{D^b(X)}({\mathscr F}^{\bullet},\,{\mathscr F}^{\bullet}).$ 
\end{Def} 
\begin{Rems}\label{R:generalization}
The inclusion functor from the category $\mathbf{Coh}(X)$ of coherent sheaves on $X$ to its derived category $D^b(X)$ is fully faithful \cite[p.\,164]{Ge-Ma}. Therefore, Definition \ref{Def-brane} coincides with Definition \ref{Def2} when it is applied to a brane consisting of  only one nontrivial term ${\mathscr F }.$    

\end{Rems}

\begin{Rems}
By the adjunction isomorphism
$${\rm Hom}_{D^b(X^{(1)})}\big(L\pi_1^*{\mathscr F}^{\bullet},\, L\pi_2^*{\mathscr F}^{\bullet}  \big)\simeq
{\rm Hom}_{D^b(X)}({\mathscr F}^{\bullet},\, {\mathscr J}^1({\mathscr F}^{\bullet}))$$  
a gauge field on ${\mathscr F}^{\bullet}$ can be considered as an element of 
 $${\rm Hom}_{D^b(X^{(1)})}\big(L\pi_1^*{\mathscr F}^{\bullet},\, L\pi_2^*{\mathscr F}^{\bullet}  \big),$$ which lifts the identity on $({\mathscr F}^{\bullet},\,\delta^{\bullet})$, in accordance with  the observation at the beginning of Subsection \ref{ss:FirstJet} (see (\ref{Deligne_Gauge})).
 
\end{Rems}

\begin{Rems}
 The gauge fields on ${\mathscr F}^{\bullet}$ are elements  
${\rm Ext}^0({\mathscr F}^{\bullet},\,{\mathscr J}^1({\mathscr F}^{\bullet}));$
i.e., the are particular open strings between ${\mathscr F}^{\bullet}$ and ${\mathscr J}^1({\mathscr F}^{\bullet})$ with 
ghost number $0$ \cite[Sect. 5.2]{Aspin}, \cite{Katz-Sharpe}.

\end{Rems}

The exact sequence (\ref{ExactJ1})
determines a distinguished triangle
$$\Omega^1({\mathscr F}^{\bullet})\overset{\iota}{\rightarrow}{\mathscr J}^1({\mathscr F}^{\bullet} ) \overset{\pi}{\rightarrow} {\mathscr F}^{\bullet}\overset{+1}{\to}$$
in the category $D^b(X)$ \cite[p.\,46]{Kas-Sch}, \cite[p.\,157]{Ge-Ma}. As   ${\rm Hom}_{D^b(X)}({\mathscr F}^{\bullet},\,.\,)$ is a cohomological functor, it follows
that 
\begin{align} 0\to&{\rm Hom}_{D^b(X)}({\mathscr F}^{\bullet},\, \Omega^1({\mathscr F}^{\bullet}))\overset{\lambda}{\longrightarrow}
 {\rm Hom}_{D^b(X)}({\mathscr F}^{\bullet},\, {\mathscr J}^1({\mathscr F}^{\bullet})) 
 \overset{\mu}{\longrightarrow} \notag \\
&{\rm Hom}_{D^b(X)}({\mathscr F}^{\bullet},\,{\mathscr F}^{\bullet})  
\overset{\nu}{\longrightarrow} {\rm Ext}^1({\mathscr F}^{\bullet},\, \Omega^1({\mathscr F}^{\bullet}))\to  \notag
\end{align}
is an exact sequence.


Since $\mu(\psi)=\pi\circ \psi$, the existence of a holomorphic gauge field $\psi$ on ${\mathscr F}^{\bullet}$, is equivalent to
$1_{{\mathscr F}^{\bullet}}\in{\rm Im}(\mu)={\rm Ker}(\nu)$. 
If $\psi$ and $\psi_1$ are gauge fields on the coherent sheaf ${\mathscr F^{\bullet}}$, then $\mu(\psi_1-\psi)=0$; i. e. $\psi_1-\psi\in {\rm Im}(\lambda)$. Hence, the set of holomorphic gauge fields on ${\mathscr F}^{\bullet}$, if nonempty, is an affine space with vector space  
 ${\rm Hom}_{D^b(X)}({\mathscr F}^{\bullet},\Omega^1({\mathscr F}^{\bullet})).$
 Thus,  the space of holomorphic gauge fields on   ${\mathscr F}^{\bullet}$, is a finite dimensional affine space.


 The Atiyah class  $a({\mathscr F}^{\bullet})$ of ${\mathscr F}^{\bullet}$  is the image of $1\in {\rm Hom}_{D^b(X)}({\mathscr F}^{\bullet},\,{\mathscr F}^{\bullet})$  in ${\rm Ext}^1({\mathscr F}^{\bullet},\, \Omega^1({\mathscr F}^{\bullet}))$. Hence, we have the following 
 proposition.
\begin{Prop}\label{P:Atiyah}
The vanishing of $a({\mathscr F}^{\bullet})$
is a necessary and sufficient condition for the existence of gauge fields on the brane ${\mathscr F}^{\bullet}.$ Furthermore, the set of gauge fields on  ${\mathscr F}^{\bullet},$ if is nonempty, is an affine space over the finite dimensional vector space ${\rm Ext}^0({\mathscr F}^{\bullet},\,\Omega^1({\mathscr F}^{\bullet})).$
\end{Prop}

 
\subsection{$B$-branes on stratified manifolds}\label{Ss:Flag}
When the manifold 
$X$ admits a certain type of stratification, the derived category $D^b(X)$
 is generated by a finite set of coherent sheaves, and the space of gauge fields on any brane over $X$ has cardinal $\leq 1.$ In this section, we will address this issue. 
  In the next   subsection, we review some properties of the generators of   $D^b(X).$ 
  Finally, we consider the  particular case where $X$
 is the variety of complete flags in ${\mathbb C}^3,$ 
  showing that the cardinality of the set of holomorphic gauge fields on any brane over $X$ is $<2.$

\subsubsection{Generators of a derived category}\label{S:Generators}

Let ${\mathbf C}$ be an abelian category.  If $A $ a complex in ${\mathbf C}$,
as usual, we denote by ${A}[l]$, with $l\in{\mathbb Z}$, the complex ${A}$ shifted $l$ to the left.  Let ${\sf E}$ be a finite set 
of objects of ${\mathbf C}$ which generates   the derived category
$D^b({\mathbf C}).$ 
 Given  $E,\tilde E$ elements of the generating set ${\sf E},$  
  let us consider morphisms $h$ between ${  E}':={  E}[l]$ and
  ${\tilde E}':={\tilde E}[l']$. We denote by ${\rm Cone}(h)={E}'[1]\oplus {\tilde E}'$ the mapping cone of $h$
   \cite[p.\,154]{Ge-Ma}. We define ${\sf E}^{(1)}$ the set obtained adding to ${\sf E}$ the elements of the form ${\rm Cone}(h)$. Hence, an element of ${\sf E}^{(1)}$ is a complex whose term at a position $p$ is either  $0$, or $E$, or a direct sum of ${E}\oplus \tilde E,$  with $E,\tilde E\in{\sf E}.$  

Repeating the process with the elements of ${\sf E}^{(1)}$ one obtains ${\sf E}^{(2)}$, etc. The objects of the triangulated subcategory generated by the family ${\sf E}$ are elements which belong to some ${\sf E}^{(m)}$.
      Therefore, an object  of the triangulated subcategory of $D^b({\mathbf C})$  generated by ${\sf E}$
 is a complex $(  G ,\,d )$, where ${ G}^{p}$ is   of the form 
 \begin{equation}\label{bigoplus}
 G^{p}=\bigoplus_{i\in I_p} E_{pi},
 \end{equation}
 with $E_{pi}\in{\sf E}$ and $i$ varying in a finite set $I_p$. (When $i$ ``runs over the empty set'', the direct sum is taken to be $0$). As ${\sf E}$  is a set generating  $D^b({\mathbf C}),$ each object of this category is isomorphic to one of the form (\ref{bigoplus}).

\smallskip

 Given two 
 complexes $(A,\,\partial_A)$
  and  $(B,\,\partial_B)$ in an additive category ${\mathbf C}$, the complex
${\rm Hom}^{\bullet}(A,\,B)$ is defined by (see \cite[p.\,17]{Iversen})
\begin{equation}\label{deltamg}
{\rm Hom}^m(A ,\,B )=\prod_{q\in{\mathbb Z}}{\rm Hom}_{\mathbf C}\big(A^q,\, B^{q+m}   \big),
\end{equation}
with the differential  $\delta_H$.
\begin{equation}\label{deltamg1}
(\delta_H^mg)^q=\partial_B^{m+q}g^q+(-1)^{m+1}g^{q+1}\partial_A^q,
\end{equation}
with $g\in {\rm Hom}^m(A ,\,B ).$ 

Denoting by ${\rm Com}({\mathbf C})$ the abelian category of complexes on ${\mathbf C}$ and fixed the object $A,$ one has the functor $F_A:{\rm Com}(\mathbf{C})\to {\rm Com}(\mathbf{Ab}),$ where $\mathbf{Ab}$ is the category of abelian groups, defined by
$$F_A(B)={\rm Hom}^{\bullet}(A,\,B),\;\; F_A(f)^m(g)=\big(f^{m+q}g^q\big),$$
$f$ a morphism from $(B,\,\partial_B)\to (C,\,\partial_C).$ Assumed that ${\mathbf C}$ has enough injectives, as $F_A$ is a left exact functor, one has the corresponding derived functor
$RF_A:D^+{(\mathbf C})\to D^+(\mathbf{Ab}).$

Since the $0$th derived functor of a left exact functor is isomorphic to the functor itself, $F_A=R^0F_A ={\rm Ext}^0(A,\,.\,).$ 
On the other hand,  ${\rm Hom}_{D({\mathbf C})}\big({A} ,\,
 B \big)={\rm Ext}^0(A,\,{  B} )$ \cite[Sect.\,10.7]{Weibel}, so
\begin{equation}\label{HomDb}
  {\rm Hom}_{D({\mathbf C})}\big({A} ,\,
 B \big)= 
{\rm Hom}^{\bullet} ({A} ,\, {B}).
\end{equation}

   One has the following lemma.
   \begin{Lem}\label{L:Hom}
   Let  ${G}$ be the complex  in the category ${\mathbf C}$ defined in (\ref{bigoplus})  and $\Hat G$ the complex ${\Hat G}^q=\bigoplus_{j\in J_q} E_{qj}.$ If ${\rm Hom}_{\mathbf C}\big( E_{pi},\, E_{qj}  \big)=0$ for all  $p,q,i,j,$ then
 $${\rm Hom}_{D^b({\mathbf C})}\big({  G} ,\, \Hat G 
  \big)=0.$$ 
  
	\end{Lem}
 {\it Proof.} From (\ref{deltamg})  together with the additivity of the functor ${\rm Hom}_{\mathbf C}(\,.\,,\,.\,),$ it follows
 $${\rm Hom}^m(G ,\,\Hat G )=\prod_p\bigoplus_{ij}{\rm Hom}_{\mathbf C}(E_{pi},\,E_{(p+m)j})=0.$$
By  (\ref{HomDb}), ${\rm Hom}_{D^b({\mathbf C})}({G} ,\,
 \Hat G )= 0.$
 \qed
 
 We shall deal with the particular case where the category ${\mathbf C}$  is a category $\mathbf{Coh}(X)$ of coherent sheaves. More precisely, let $X$ be a smooth complex $n$-dimensional  variety, let 
  $\mathbf{Qcoh}(X)$ denote the category of quasi-coherent sheaves on $X.$ One has the derived functor 
$ RF_A:D^+(\mathbf{Qcoh}(X))\to D^+(\mathbf{Ab}).$ 
Since ${\rm Ext}^i(A,\, B)=0$ for $A,B$ complexes of coherent sheaves and 
$i\notin [0,\,n],$ we have the functor
$$RF_A:D^b(X):=D^b(\mathbf{Coh}(X))\to D^b(\mathbf{Ab}).$$

From Lemma \ref{L:Hom} applied to the category $D^b(X)$, it follows the following corollary.
 
 \begin{Cor}\label{Coro:1}
Let ${\mathscr G}^p=\bigoplus_{pi}{\mathscr E}_{pi}$  be an object of $D^b(X),$
 where ${\mathscr E}_{pi}\in\mathbf{Coh}(X),$ and $\Hat{\mathscr  G}$ is an object of $D^b(X)$
 category isomorphic the complex $ \Omega^1_X \otimes_{{\mathscr O}_X}{\mathscr G}^p,$  then 
   ${\rm Hom}_{D^b(X)}\big({ {\mathscr  G}} ,\, \Hat{\mathscr  G }
  \big)=0,$
  if ${\rm Hom} ({\mathscr E}_{ip},\,\Omega^1({\mathscr E}_{jq})) =0$ for all $p,q,i,j.$ 
  \end{Cor}

  
\smallskip

\paragraph{{\bf Generators of the derived category $D^b(X)$}.}

Let $X$ be a locally Noetherian scheme, $i: Z\hookrightarrow X$ a closed subscheme of $X$ and $j:U\hookrightarrow X$ an open part of $X.$ Since $j$ is an open embedding, the functor proper inverse image $j^!$ is defined between the categories of modules $\mathbf{Mod}({\mathscr O}_X)\to\mathbf{Mod}({\mathscr O}_U).$ 
 One has the pair of adjoint functors $(j_!\dashv j^!=j^*)$ and the pair $(i^*\dashv i_*=i_!).$ 
 $$ \mathbf{Mod}({\mathscr O}_X)\overset{j^!}{\to}\mathbf{Mod}({\mathscr O}_U)
 \overset{j_!}{\longrightarrow}\mathbf{Mod}({\mathscr O}_X)
   \overset{i^*}{\to}\mathbf{Mod}({\mathscr O}_Z)
 \overset{i_*}{\longrightarrow}\mathbf{Mod}({\mathscr O}_X).$$
 

For every  sheaf ${\mathscr F}$ on $X,$ the corresponding adjunction morphisms give rise to the sequence of morphisms of ${\mathscr O}_X$-modules
$\,j_!j^!{\mathscr F}\to{\mathscr  F}\to i_*i^*{\mathscr F}.$
 If ${\mathscr F}$ is a coherent ${\mathscr O}_X$-module, then sheaf $i^*{\mathscr F}$ is coherent. Since $i$ is a proper map, $i_*i^*{\mathscr F}$ is a coherent ${\mathscr O}_X$-module.  
As $j_!j^!{\mathscr F}$ is a subsheaf of ${\mathscr F}$, it is also coherent.
 
 In particular, when $U=X\setminus Z$, for any ${\mathscr F}\in\mathbf{Coh}(X)$, one has the exact sequence of ${\mathscr O}_X$-modules
$ 0\to j_!j^!{\mathscr F}\to{\mathscr  F}\to i_*i^*{\mathscr F}\to 0.$
  If ${\mathscr F}^{\bullet}$ is an object of the derived category $D^b(X),$  then we have
the distinguished triangle 
$Rj_!j^!{\mathscr F}^{\bullet}\to{\mathscr  F}^{\bullet}\to Ri_*i^*{\mathscr F}^{\bullet}\overset{+1}{\longrightarrow},$
where $j^!{\mathscr F}^{\bullet}\in D^b(\mathbf{Coh}(U))$ and $i^*{\mathscr F}^{\bullet}\in D^b (\mathbf {Coh}(Z)).$ Hence, as $j_!$ and $i_*$ are exact  functors, one has the
 following proposition.
\begin{Prop}\label{P:gk} 
If $\{g_k\}$ (resp. $\{g'_{k'}\}$) is a set 
  generators of the derived category $D^b(\mathbf{Coh}(U))$ (resp. $D^b(\mathbf {Coh}(Z))$) as a triangulated category, then $\{j_!g_k,\,i_*g'_{k'}\}$ is a set that generates  $D^b(X).$
  \end{Prop}

If  $U={\rm Spec}\,R$, where $R$ is  a Noetherian ring, the  category $\mathbf{Coh}(U)$ of coherent sheaves on $U$ is equivalent to the category of finitely generated $R$-modules, with  $R=\Gamma(U,\,{\mathscr O}_U).$   The functor $\Gamma(U,\,\cdot\,)$ defines the equivalence.

Moreover, if $R$ is a regular local ring, then its global dimension is finite, and any   object $N$ of the category $\mathbf{Mod}_{f}(R)$  (of the  $R$-modules of   finite type) admits a finite free resolution
\begin{equation}\label{longresolution}
0\to S_m\overset{f_m}{\longrightarrow} S_{m-1} \to\dots\to S_1\overset{f_1}{\longrightarrow} S_0\to N\to 0,
\end{equation}
with $S_a$ a finite direct sum $\oplus R.$
 
Denoting by $\mathbf{T}$ the triangulated subcategory of $D^b(\mathbf{Mod}_{f}(R))$ generated by $R,$ then the $S_a$ are objects of $\mathbf{T}.$ On the other hand, if (\ref{longresolution}) reduces to
$0\to S_1\overset{f_1}{\longrightarrow} S_0\to N\to 0$ in the category $\mathbf{Mod}_{f}(R),$ then one has the 
  distinguished triangle 
$$S_1\overset{f_1}{\longrightarrow} S_0\to{\rm Cone}(f_1)\overset{+1}{\longrightarrow}.$$
Thus, $N$ belongs to the category $\mathbf{T}.$ In the general case where case   $N$ has a resolution as (\ref{longresolution}), then one can consider the exact sequences
$$0\to S_m\overset{f_m}{\longrightarrow} S_{m-1} \to {\rm Im}(f_{m-1})\to 0,$$ 
and
$$ 0 \to {\rm Im}(f_{m-1})\to S_{m-2}\to\dots \to N\to 0.$$
By the previous result, ${\rm Im}(f_{m-1})$ belongs to $\mathbf{T}.$ By induction, it follows that $N$ is also an object of $\mathbf{T}.$ Hence
$\mathbf{Mod}_f(R)\subset{\mathbf T}.$ As $D^b(\mathbf{Mod}_f(R))$ is generated by $\mathbf{Mod}_f(R),$ it follows that $R$ generates the derived  category $ D^b(\mathbf{Mod}_f(R)).$ 

From the equivalence between $\mathbf{Coh}(U)$ and $\mathbf{Mod}_f(R),$ one deduces   the following proposition.
\begin{Prop}\label{P:AffineScheme}
If $U$ is a Noetherian affine scheme, then 
  ${\mathscr O}_U$ generates $D^b(U)$, the bounded derived category of coherent sheaves on $U$.
  \end{Prop}
  From Proposition \ref{P:gk}, it follows the following  corollary (see \cite[Lem. 2.3.2]{BOH}). 
\begin{Cor}\label{j!O_X}
If $U\overset{j}{\hookrightarrow}X$ is an open part of $X$  isomorphic to a Noetherian affine scheme and $Z=X\setminus U,$ then the derived category $D^b(X)$ is generated by $j_!{\mathscr O}_U$ together with  the elements of   
$i_*\{{\rm generators\; of}\, D^b(Z) \},$
where $i$ is the inclusion $Z\hookrightarrow X.$
\end{Cor}


\smallskip


 Let $X$ be a smooth a complex $n$-dimensional variety. Let $Z_1$ be  an irreducible  subvariety of $X$ of codimension $1,$ and $i_1:Z_1\hookrightarrow X$ the corresponding inclusion. We set
$U_1:=X\setminus  Z_1$ and  $j_1:U_1\hookrightarrow X.$  
 Let us assume that $U_1$ is isomorphic to an affine variety. 
The exact sequence 
$$0\to j_{1!}j_1^{!}{\mathscr O}_X\to{\mathscr O}_X\to i_{1*}i^*_1{\mathscr O}_X \to 0,$$
can be written
$0\to {\mathscr M}_1\to{\mathscr O}_X\to {\mathscr O}_{Z_1}\to 0,$ where ${\mathscr M}_1$ is the  invertible  sheaf ${\mathscr O}_X(-Z_1),$ which is supported in $U_1,$ i.e. the stalk $({\mathscr M}_1)_x=0$ for all $x\in Z_1.$

As $j_1^!{\mathscr O}_X={\mathscr O}_{U_1}$, then ${\mathscr M}_1=j_{1!}{\mathscr O}_{U_1}$. From Corollary \ref{j!O_X}, it follows that the derived category $D^b(X)$ is generated by ${\mathscr M}_1$ together with the elements of
$i_{1*}\{{\rm generators\; of}\, D^b(Z_1) \}.$

If $Z_2\overset{i_2}{\hookrightarrow} Z_1$ is an irreducible  hypersurface of $Z_1$  and $U_2:= Z_1\setminus  Z_2$ is isomorphic to an affine variety, then we have an exact sequence
$$0\to{\mathscr  M}_2\to {\mathscr O}_{Z_1}\to i_{2*}i_2^{*}{\mathscr O}_{Z_1}\to 0,$$
where ${\mathscr M}_2$ is  the   sheaf on $Z_1,$ ${\mathscr O}_{{  Z}_1}(-Z_2).$
Hence, $D^b(Z_1)$ is generated by ${\mathscr M}_2$ together with
$i_{2*} \{{\rm generators\, of}\, D^b(Z_2)  \}.$
Consequently, $D^b(X)$ is generated by
$$\{{\mathscr M}_1,\,i_{1*}{\mathscr M}_2\}\cup i_{1*}i_{2*}\big(\{{\rm generators\, of}\, D^b(Z_2)  \}\big).$$

Let us assume that there exists a tower $X=Z_0\supset Z_1\supset Z_2\supset\dots\supset Z_n,$ where $Z_k\overset{i_k}{\hookrightarrow} Z_{k-1}$ is an irreducible hypersurface   of $Z_{k-1}$ with $i=1,\dots, n.$ 
 According to  the above argument, ${\mathscr M}_k={\mathscr O}_{Z_{k-1}}(-Z_k)$ is a sheaf supported on $U_k:=Z_{k-1}\setminus Z_k$. We set $\iota_k=i_1\circ\dots\circ i_k:Z_k\hookrightarrow X$  and $\iota_0:={\rm id}:Z_0\to X.$ For $k=1,\dots, n$ we denote
${\mathscr L}_k:=\iota_{k-1*}({\mathscr M}_k).$ 
We set ${\mathscr L}_{n+1}$ for denoting the sheaf on $X$ direct sum of the skyscraper sheaves on the points of the finite set $Z_n.$ Then 
${\mathscr L}_{n+1}=\iota_{n*}\big(\{{\rm generators\, of}\, D^b(Z_n)\}).$
We have proved the following proposition.

\begin{Prop}\label{generadores DX}
 With the above notations,
the derived category
 $D^b(X)$ is generated by the set of ${\mathscr O}_X$-modules
$\,{\sf E}=\{{\mathscr L}_1,\dots, {\mathscr L}_n,\,{\mathscr L}_{n+1}\}.$
\end{Prop}

\begin{Rems}\label{jandk}
 When $j\ne k$ the intersection $U_j\cap U_k=\emptyset.$ Hence, ${\mathscr L}_j$ and $\Omega^1({\mathscr L}_k)$ have disjoint supports; thus,   
\begin{equation}\label{LjLk}
{\rm Hom} \big({\mathscr L}_j,\,\Omega^1_X\otimes_{{\mathscr O}_X}{\mathscr L}_k\big)=0.
\end{equation}
\end{Rems}

\smallskip

\begin{Lem}\label{Lemitems}
With the notation introduced above,
$ {\rm Hom}({\mathscr L}_k,\,{\mathscr L}_k\otimes_{{\mathscr O}_X}\Omega^1_X) =H^{1,0}(Z_k),$ for $k=1,\dots,n.$
\end{Lem}
{\it Proof.}
 We set in this proof  ${\mathscr L}:={\mathscr L}_k,$ ${\mathscr M}:={\mathscr M}_k,$ and $\iota:=\iota_{k-1}$  for the inclusion
$Z:=Z_{k-1}\hookrightarrow X;$ thus ${\mathscr L}=\iota_*{\mathscr M}.$ We will also delete  the subscript $X$ in the sheaves  ${\mathscr O}_X=:{\mathscr O}$ and
$\Omega^1_X=:\Omega^1.$ 
By ${\mathscr Hom}\,_{\mathscr R}(\,.\, ,\,.\,)$ we denote the sheaf of linear homomorphisms between the  corresponding ${\mathscr R}$-modules.

By \cite[Prop.\,2.5.13]{Kas-Sch} and \cite[Cor.\,2.3.4]{Kas-Sch} 
$$ {\mathscr Hom}\,_{\mathscr O}({\mathscr L},\,{\mathscr L}\otimes_{\mathscr O}\Omega^1)
\simeq\iota_* {\mathscr Hom}\,_{{\mathscr O}_Z}({\mathscr M},\,{\mathscr M}\otimes_{{\mathscr O}_Z}\Omega^1_Z).$$

By \cite[Prop.\,2.5.13]{Kas-Sch} and \cite[Cor.\,2.3.4]{Kas-Sch} 
$$ {\mathscr Hom}\,_{\mathscr O}({\mathscr L},\,{\mathscr L}\otimes_{\mathscr O}\Omega^1)
\simeq\iota_* {\mathscr Hom}\,_{{\mathscr O}_Z}({\mathscr M},\,{\mathscr M}\otimes_{{\mathscr O}_Z}\Omega^1_Z).$$
Hence, as ${\mathscr M}$ is an invertible ${\mathscr O}_Z$-module,
\begin{align}
\notag 
&{\rm Hom}({\mathscr L},\,{\mathscr L}\otimes_{\mathscr O}\Omega^1)\simeq
\Gamma\big(X,\,\iota_*{\mathscr Hom}\,_{{\mathscr O}_Z}({\mathscr M},\,{\mathscr M}\otimes_{{\mathscr O}_Z}\Omega^1_Z)\big)= \\ \notag
&\Gamma\big(Z,\,{\mathscr Hom}\,_{{\mathscr O}_Z}({\mathscr M},\,{\mathscr M}\otimes_{{\mathscr O}_Z}\Omega^1_Z)\big)\simeq {\rm Hom}_{{\mathscr O}_Z}({\mathscr O}_Z,\Omega^1_Z)=H^{1,0}(Z).\notag
\end{align} \qed

\begin{Thm}\label{Th:Z_i} Let $X$ be an $n$-dimensional  smooth variety and suppose that: 
\begin{enumerate}
\item There exists $X=Z_0 \supset Z_1\supset\dots\supset Z_n,$ where $Z_{i+1}$ is a smooth irreducible subvariety of $Z_{i},$  with ${\rm codim}_{Z_i}Z_{i+1}=1$. 
\item $U_j=Z_{j-1}\setminus Z_{i}$ is isomorphic to an affine variety.
\end{enumerate}
If $H^{1,0}(Z_j)=0$ for all $j,$ then the set of holomorphic gauge fields on any $B$-brane over $X$ has cardinal $\leq 1.$
\end{Thm}
{\it Proof.} From Lemma \ref{Lemitems} together with Remark \ref{jandk}, it follows that ${\rm Hom}({\mathscr L}_j,\,\Omega^1({\mathscr L}_k))=0,$ for all $j,k.$ Given a brane ${\mathscr F}^{\bullet}$ over $X,$ from Corollary  \ref{Coro:1} and Proposition \ref{generadores DX}, we deduce that ${\rm Hom}_{D^b(X)}({\mathscr F}^{\bullet},\,\Omega^1({\mathscr F}^{\bullet}))=0.$ Now the theorem is consequence of   Proposition \ref{P:Atiyah}.
 \qed

 Using this theorem when $X$
 is the projective space,  it has been proven in \cite{Vina24} that
  the set of holomorphic gauge fields on any $B$-brane on ${\mathbb P}^n$ has cardinal $\leq 1.$




\subsubsection{Flag varieties}\label{S:Flag_varieties}   

In this section, we will first recall some properties of the
  Schubert varieties. We denote by $X:={\rm GL}(n,\,{\mathbb C})/B$, where $B$ is the Borel subgroup of ${\rm GL}(n,\,{\mathbb C})$ consisting of the upper triangular matrices. $X$ is the manifold of complete flags in ${\mathbb C}^n.$
  
   Let $T$ denote the maximal torus of
 ${\rm GL}(n,\,{\mathbb C})$ formed by the
 diagonal matrices.
    We set $W:=N(T)/T$, where $N(T)$ is the normalizer of $T$ in the group ${\rm GL}(n,\,{\mathbb C}).$ As it is well-known, the Weyl group $W$ is isomorphic to the symmetric group $S_n$. The Bruhat   partial order \cite[Sect.\,10.5]{Fulton-Y} in $S_n$ will be denoted by $\leq$. 
 
   Given $w\in W,$ we set $C^0_w$  for the Schubert cell   of $X$ associated to $w$ \cite{Fulton-Y}. The closure (in Zariski topology) of this cell is the corresponding Schubert subvariety of $X$
	\begin{equation}\label{S_w}
	C_w:=\bigcup_{v\leq w}C_v^0.
	\end{equation}
	
   The Poincar\'e polynomial of $C_w$ given in \cite {A-C} shows that the Betti number $b_1(C_w)$
    vanishes, if $C_w$ is any smooth Schubert variety.
  
   Theorem 1 of \cite{L-S} gives a criterion to determine whether a Schubert variety in a flag variety is singular.  Let $w$ be the permutation $w=(t_1,\dots,t_n)$. The variety $C_w$ is singular iff there are integers $i,j,k,l$ with $1\leq i<j<k<l\leq n$ such that, 
      either
$$t_k<t_l<t_i<t_j \quad \text{or} \quad    t_l<t_j<t_k<t_i.$$        
  In particular, all the Schubert varieties  of $X={\rm GL}(3,\,{\mathbb C})/B$ are smooth. Thus,  if $C_w$ is such a variety, as  the $b_1(C_w)=0$, the Hodge numbers $h^{1,0}(C_w)$ and $h^{0,1}(C_w)$ vanish
  \begin{equation}\label{b1h10}
  h^{1,0}(C_w)=h^{0,1}(C_w)=0.
  \end{equation}
   
   From now on in this section, $X$   will denote the variety of complete flags in ${\mathbb C}^3.$  The Bruhat order in the group $S_3$ is showed in the following diagram \cite{Brion}.
	
	
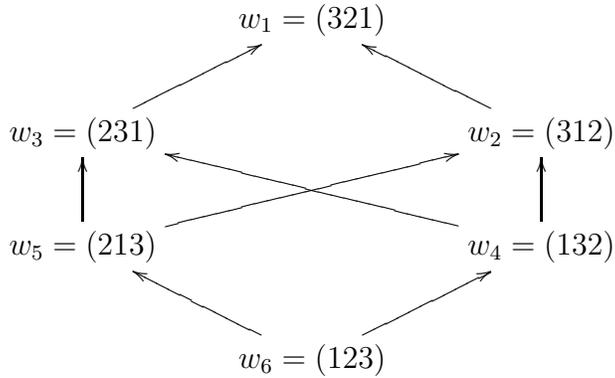
	$$\xymatrix{
{} & w_1=(321)  & {} \\   
w_3=(231)\ar[ur] & {} &  w_2=(312) \ar[ul]    \\   
w_5=(213)\ar[u] \ar[urr] & {} & w_4=(132) \ar[u] \ar[ull] \\
{} & w_6=(123)\ar[ul] \ar[ur]  & {} \\
}
  $$


	Hence,
	there are only two pairs of permutations that are not comparable in the Bruhat order, namely, $(w_2,\, w_3)$ and $(w_4,\, w_5).$ For any other pair $(w,\,w')$ the intersection $C_w\cap C_{w'}=C_v$, where $v={\rm min}\{w,\,w'\}.$

The length $l(w)$ of the permutation $w$ is defined by
	$$l(w)=\#\{ a<b\,|\, w(a)>w(b)\}.$$
	We will set $C^0_k$ for the cell $C^0_{w_k}$ and $C_k:=C_{w_k}.$
	The cell $C^0_{w}$ is an affine variety isomorphic to ${\mathbb C}^{l(\omega)}$ \cite[p.\,157]{Fulton-Y}. Hence, 
	${\rm dim}\, C_4=1={\rm dim}\, C_5$ and ${\rm dim}\, C_2=2={\rm dim}\, C_3.$ Furthermore, $X$ is the disjoint union of the cells
	$$X=\bigsqcup_{k=1}^6C^0_k.$$
	On the other hand, the inclusion relations among the Schubert subvarieties of $X$ can be obtained from (\ref{S_w}) by means of the Bruhat order   depicted above.
	

We consider the following tower of subspaces of $X$
$$X=Z_0=C_1\overset{i_1}{\hookleftarrow} Z_1=C_2\cup C_3  \overset{i_2}{\hookleftarrow} Z_2=C_4\cup C_5 \overset{i_3}{\hookleftarrow}Z_3=C_6$$
That is, $Z_i$ is the union of the Schubert varieties of dimension $3-i.$ We set $U_i:=Z_{i-1}\setminus Z_i$.  From (\ref{S_w}), it follows the expression of each $U_i$ in terms of Schubert cells. 
One has the the following inclusion relations
$$Z_0\overset{j_1}{\hookleftarrow} U_1=C^0_1,\;\;\; Z_1\overset{j_2}{\hookleftarrow} U_2=C^0_2\cup C^0_3,\;\;\; Z_2\overset{j_3}{\hookleftarrow} U_3=C^0_4\cup C^0_5.   $$

A family of generators for $D^b(X)$ will then be constructed by successive applications of   Corollary \ref{j!O_X}.

\smallskip

\noindent 
{\sc $U_1$ open part of $Z_0=X.$} In the flag variety $X$, we have the divisor $Z_1=X\setminus C^0_1=C_2\cup C_3,$ and the  partition of $X$ defined by the affine variety $U_1=C^0_1$ and $Z_1.$ We denote by $j_1$ and $i_1$ the respective inclusions 
$$U_1\overset{j_1}{\hookrightarrow} X\overset{i_1}{\hookleftarrow} Z_1.$$
 On the other hand, we have the exact sequence of ${\mathscr O}_X$-modules
$$0\to{\mathscr  M}_1:=j_{1!}j^{!}_1{\mathscr O}_X\to{\mathscr O}_X\to i_{1*}i^*_1{\mathscr O}_X\to 0.$$
The nontrivial term on the left satisfies
\begin{enumerate}
\item ${\mathscr M}_1$ is the invertible sheaf ${\mathscr O}_X(-C_2-C_3)$ associated to the divisor $-C_2-C_3.$
\item ${\rm Supp}\, {\mathscr M}_1\subset U_1.$
\item  As $j^!_1{\mathscr O}_X={\mathscr O}_{U_1}$, from Corollary \ref{j!O_X}, it follows that 
 $D^b(X)$ is generated by ${\mathscr M}_1$ together with  
$i_{1*}\{{\rm generators \; of}D^b(Z_1)\}.$ We write
\begin{equation}\label{Db(X)=}
D^b(X)=\langle {\mathscr M}_1,\,i_{1*} D^b(Z_1)   \rangle.
\end{equation}
\end{enumerate}

\smallskip

\noindent 
{\sc $U_2$ open subset of  $Z_1.$} Next we will consider a divisor $Z_2=C_4+C_5$ of $Z_1$	 and we apply to it the above analysis. 
We will first state the following simple lemma.

\begin{Lem}\label{L:Flag}
Let $Y$ be a complex variety, and $Y_1,\,Y_2$ disjoint affine subvarieties of $Y,$ 
 then the derived category $D^b(Y)$ is generated by 
$\{\alpha_{1!}{\mathscr O}_{Y_1},\, \alpha_{2!}{\mathscr O}_{Y_2}\},$
where 
$\alpha_k:Y_k\hookrightarrow Y$  denotes the corresponding embedding.
\end{Lem}
{\it Proof.}
Given ${\mathscr F}\in\mathbf{Coh}(Y),$ one has the coherent sheaves
 ${\mathscr F}_k:=\alpha_{k!}\alpha_k^{!}{\mathscr F}$ and 
 the exact sequence of ${\mathscr O}_Y$-modules \cite[Prop.\,2.3.6 (vii)]{Kas-Sch}
 $$0\to{\mathscr  F}_{Y_1\cap Y_2}=0\to{\mathscr  F}_1\oplus {\mathscr F}_2\to {\mathscr F}\to 0.$$
 By 
 Proposition \ref{P:AffineScheme},
 $D(Y_k)$ is generated  by 
${\mathscr O}_{Y_k}.$ From the above exact sequence, it follows that the sheaf ${\mathscr F}$ belongs to the subcategory of $D^b(Y)$ generated by $\alpha_{1!}{\mathscr O}_{Y_1},\,$  $\alpha_{2!}{\mathscr O}_{Y_2}.$ 
\qed

\smallskip

Going back to  the variety $X$.
Next we will consider a divisor $Z_2=C_4+C_5$ of $Z_1.$ Then $U_2:=Z_1\setminus Z_2=C_2^0\cup C_3^0.$ 	For $k=2,3$ we have the inclusions
$$C^0_k\overset{\alpha_k}{\hookrightarrow} U_2
\overset{j_2}{\hookrightarrow} Z_1\overset{i_2}{\hookleftarrow}Z_2.$$
 As $C^0_2\cap C^0_3 =\emptyset,$ 
from Lemma \ref{L:Flag} together with   Corollary \ref{j!O_X}, it follows that    the set $\{ \alpha_{k!}{\mathscr O}_{C^0_k}\}_{k=2,3}$ is a set of generators of $D^b(U_2).$

For $k=2,3$ the ${\mathscr O}_{Z_1}$-module
 $(j_2\alpha_k)_!{\mathscr O}_{C^0_k},$
 will be denoted   ${\mathscr M}_k.$ 
\begin{enumerate}
\item Since $Z_1\setminus C^0_k=C_{k'},$ with $k'\in\{2,\,3\}\setminus \{k\},$  then ${\mathscr M}_k$ is the invertible sheaf in $Z_1$ associated to the divisor $-C_{k'}.$  That is, ${\mathscr M}_k={\mathscr O}_{Z_1}(-C_{k'}).$
\item ${\rm Supp}\,{\mathscr  M}_k\subset C^0_k .$
\item From 
 Corollary \ref{j!O_X}, we conclude
 $$D^b(Z_1)=\langle {\mathscr M}_2,\, {\mathscr M}_3,\, i_{2*}D^b(Z_2)\rangle.$$
\end{enumerate}

 From (\ref{Db(X)=}), it follows 
 \begin{equation}\label{Db(X)2}
 D^b(X)=\langle{\mathscr  M}_1,\, i_{1*}{\mathscr M}_2,\, i_{1*}{\mathscr M}_3,\, i_{1*}i_{2*}D^b(Z_2)   \rangle.
 \end{equation}
 
 \smallskip
 
\noindent 
{\sc $U_3$ open subset of  $Z_2.$}
In a similar way, we define $U_2=C^0_4\cup C^0_5$ and $Z_3=Z_2\setminus U_2=C_6.$ For $k=4,5$ one has the embeddings 
 $$C^0_k\overset{\alpha_k}{\hookrightarrow}U_3\overset{j_3}{\hookrightarrow} Z_2\overset{i_3}{\hookleftarrow}Z_3.$$ 
The ${\mathscr O}_{Z_2}$-module
 $(j_3\alpha_k)_!{\mathscr O}_{C^0_k},$
 will be denoted ${\mathscr M}_k.$ 
\begin{enumerate}
\item Since $Z_2\setminus C^0_k=C_{k'},$ with $k'\in\{4,\,5\}\setminus \{k\},$  then ${\mathscr M}_k$ is the invertible sheaf in $Z_2$ associated to the divisor $-C_{k'}.$  That is, ${\mathscr M_k}={\mathscr O}_{Z_2}(-C_{k'}).$
\item ${\rm Supp}\,{\mathscr  M}_k\subset C^0_k .$
\item By Corollary \ref{j!O_X},
 $D^b(Z_2)=\langle{\mathscr  M}_3,\, {\mathscr M}_4,\, i_{3*}D^b(Z_3)\rangle.$
\end{enumerate}
From (\ref{Db(X)2}), it follows 
 \begin{equation}\label{Db(X)3}
 D^b(X)=\langle {\mathscr M}_1,\, i_{1*}{\mathscr M}_2,\, i_{1*}{\mathscr M}_3,\, (i_1i_2)_*{\mathscr M}_4,\,(i_1i_2)_*{\mathscr M}_5,\,(i_1i_2i_3)_*D^b(Z_3)   \rangle.
 \end{equation}

For the sake of simplicity in the notation we set
$$ {\mathscr L}_1:= {\mathscr M}_1,\,{\mathscr L}_2:= i_{1*}{\mathscr M}_2,\, {\mathscr L}_3:=i_{1*}{\mathscr M}_3,\, {\mathscr L}_4:=(i_1i_2)_*{\mathscr M}_4,\,{\mathscr L}_5:=(i_1i_2)_*{\mathscr M}_5.$$

\smallskip

\noindent
{\sc The point $Z_3.$}
 ${\mathscr L}_6$ will denote the skyscraper sheaf on $X$ at the point $C_6.$ Thus, 
 \begin{Prop}\label{P:D(X)final}
 The derived category $D^b(X)$ is generated by the set $\{{\mathscr L}_1,\dots, {\mathscr L}_6\}.$
 \end{Prop}

Since ${\rm Supp}({\mathscr L}_k)\subset C^0_k$ and 
 $C^0_r\cap C^0_s=\emptyset.$ for   $r\ne s.$  Then 
 \begin{equation}\label{LrLs}
 {\rm Hom}_{{\mathscr O}_X}\big({\mathscr L}_r,\,{\mathscr L}_s\otimes_{{\mathscr O}_X}\Omega^1_X\big)=0.
 \end{equation}


\begin{Prop}\label{Lr-LrO}
${\rm Hom}_{{\mathscr O}_X}({\mathscr L}_r,\,{\mathscr L}_r\otimes_{{\mathscr O}_X}\Omega^1_X)=0,\,$ for $r=1,\dots,6.$
\end{Prop}
{\it Proof.} As in the proof of Lemma \ref{Lemitems},
${\rm Hom}_{{\mathscr O}_X}({\mathscr L}_1,\,{\mathscr L}_1\otimes_{{\mathscr O}_X}\Omega^1_X)=
H^{1,0}(Z_1).$  
 On the other hand, $Z_1=C_2\cup C_3$,  and by (\ref{b1h10}) $H^{1,0}(C_2)=H^{1,0}(C_3)=0.$ Thus, 
 $H^0(Z_1,\,\Omega^1_{Z_1})=0.$  That is,  the proposition is proved for the case $r=1. $ The  proof of the remaining cases is analogous.
\qed


\smallskip
\begin{Thm}\label{Thm:flag}
The set of holomorphic gauge fields on any $B$-brane over the variety of complete flags in ${\mathbb C}^3$ is $\leq 1$.
\end{Thm}
{\it Proof.}    
An object of $D^b(X)$ is isomorphic to one of the   triangulated category generated by the family ${\mathscr L}_1,\dots, {\mathscr L}_6$ in $D^b(X).$ 
 Let $({  G}^{\bullet},\,d^{\bullet})$ be an object of this triangulated category.
 Then ${G}^p$ is a sheaf of the form $G^p=\oplus_{i} {\mathscr L}_{pi},$ with $pi\in\{1,\dots,6\}.$
 From Corollary \ref{Coro:1} together with (\ref{LrLs})  and Proposition \ref{Lr-LrO}, it follows
${\rm Hom}_{D^b(X)}\big({G}^{\bullet},\,\Omega^1(  G^{\bullet})    \big)=0.$
The theorem is a consequence of Proposition \ref{P:Atiyah}.
\qed



\section{Holomorphic Yang-Mills fields}\label{ss:holomorphicYang}

\subsection{Hermitian sheaves}\label{Ss:Hermitian}
The fiber at $x\in X$ of a coherent ${\mathscr O}_X$-module ${\mathscr G}$ 
will be denoted by 
${\mathscr G}_{(x)}:= {\mathscr G}_{x}/{\mathfrak m}_x {\mathscr G}_{x}$, where ${\mathfrak m}_x$ is the maximal ideal of ${\mathscr O}_x$.  If ${\mathcal Z}$ is a section of ${\mathscr G}$ the corresponding vector in 
	${\mathscr G}_{(x)}$ is denoted  by $Z(x)$.  
	On the other hand, the singular set ${\mathcal S}$ of ${\mathscr G}$ is a closed analytic subset of $X$ whose codimension is greater than or equal to $1$ \cite[Chap.\,V, Thm.\,5.8]{Koba}. 
Moreover,  ${\mathscr G}$ is locally free on $X\setminus{\mathcal S}$. 
We set ${ G}$ for the vector bundle over $X\setminus {\mathcal S}$, with fibers $ G(x):={\mathscr G}_{(x)}$, determined by the locally free sheaf ${\mathscr G}|_{X\setminus {\mathcal S}}$.
	
The following definition generalizes the given in \cite[Chap.\,III, Sect.\,1]{Wells} for locally free sheaves.	
 \begin{Def}\label{D:HermitianSheaf}
A Hermitian  metric on the coherent sheaf ${\mathscr G}$ is a set $\{\langle\,,\,\rangle_x\}_{x\in X}$ of
Hermitian metrics on the fibers of ${\mathscr G}$, satisfying the following condition:
Given ${\mathcal Z}_1,{\mathcal Z}_2$ two   sections of ${\mathscr G}$ on an open $U$ of $X,$  
 the map 
$${\mathcal Z}_1\cdot {\mathcal Z}_2:x\in U\mapsto \langle Z_1(x),\,Z_2(x)\rangle_{x}\in{\mathbb C}$$
 is bounded, and its restriction 
to   $U\setminus {\mathcal S}$ is
 $C^{\infty}.$ 
 A   sheaf endowed with a Hermitian metric is called   a Hermitian sheaf.
\end{Def}

	Given ${\mathcal Z}_1,{\mathcal Z}_2\in\Gamma(X,\,{\mathscr G}),\,$   ${\mathcal Z}_1\cdot {\mathcal Z}_2$ is a bounded map with possible discontinuities in the singular locus. 	If $X$ is a K\"ahler manifold, the function  ${\mathcal Z}_1\cdot {\mathcal Z}_2$ may be non-integrable on $X$ with respect to the volume form defined by the K\"ahler metric.

\smallskip 

If ${\mathscr G}$ is a {\it torsion-free} coherent sheaf on the K\"ahler manifold $X,$ then ${\rm codim}\,{\mathcal S}\geq 2$ \cite[Cap.\,V, Cor.\,5.15]{Koba}. Hence, there exists an analytic closed    subset ${\mathcal W}$ of $X,$
containing ${\mathcal S},$  and such that $1\leq{\rm codim}\,{\mathcal W}<{\rm codim}\,{\mathcal S}.$ Let ${\mathcal N}$ be a tubular neighborhood of ${\mathcal S}$ in ${\mathcal W}.$ Then
\begin{enumerate}
\item ${\mathcal Z}_1\cdot {\mathcal Z}_2$ is a continuous map on the closed space
$X\setminus{\mathcal N},$ which is disjoint with the singularity set ${\mathcal S}.$ 
\item ${\mathcal W},$  as a subset of $X$ with codimension $\geq 1,$  has measure zero with respect to the differential form ${\rm d}{\rm vol}.$
\end{enumerate}
Thus, the following integral is well-defined
\begin{equation}\label{(mathcalZ,Z)torsion}
 ({\mathcal Z}_1,\,{\mathcal Z}_2): =
 \int_{X\setminus {\mathcal N}} {\mathcal  Z}_1\cdot {\mathcal Z}_2 \,{\rm d}{\rm vol}.
 \end{equation}
From the   boundedness of ${\mathcal  Z}_1\cdot {\mathcal Z}_2$ together with  
 $(2),$  it follows that this definition does not depend on the choices of the analytic closed set ${\mathcal W}$ and the neighborhood ${\mathcal N}.$

\begin{Rems}
A Hermitian structure $\langle\,,\,\rangle$ on 
	a coherent sheaf ${\mathscr F}$ determines a Hermitian metric on the coherent sheaf
	 ${\mathscr End}({\mathscr F})$ in a natural way. On the other hand, if $X$ is a K\"ahler manifold, the metric  defines a Hermitian structure 
	 on the locally free sheaf $\Omega^k$ of holomorphic $k$-forms.
	  Thus, if ${\mathscr F}$ is a Hermitian sheaf over a K\"ahler manifold,  we have a Hermitian metric on the sheaf $\Omega^k\otimes_{{\mathscr O}_X}{\mathscr End}({\mathscr F}),$  which will also denoted
$\langle\, ,\,\rangle$.
\end{Rems} 
	
\subsubsection{Yang-Mills functional}		 
	  
Given $\nabla$ a holomorphic connection on the
 coherent sheaf $\mathscr F.$  Then the corresponding operator $\nabla$
 defines a morphism of ${\mathbb C}_X$-modules $\nabla^{(k)}:\Omega^k({\mathscr F}) {\to} \Omega^{k+1}({\mathscr F})$ in the usual way.  
The composition ${\mathcal K}_{\nabla}:=\nabla^{(1)}\circ\nabla:{\mathscr F}\to\Omega^2({\mathscr F})$ is the curvature of $\nabla$; and as it is well-known 
\begin{equation}\label{mathcalKnabla}
{\mathcal K}_{\nabla}\in{\rm Hom}({\mathscr F},\,\Omega^2({\mathscr F}))=\Gamma(X,\,{\mathscr Hom}_{{\mathscr O}_X}({\mathscr F},\, \Omega^2({\mathscr F}))).
\end{equation}
The connection is said to be {\it flat} if ${\mathcal K}_{\nabla}=0$.

	For each $x\in X,$ we denote by $\alpha_x$ and $\lambda_x$ the natural morphisms 
	$$\big({\mathscr Hom}_{{\mathscr O}_X}({\mathscr F},\, \Omega^k({\mathscr F}))\big)_x \overset{\alpha_x}{\to}
	{\rm Hom}_{{\mathscr O}_x}\big({\mathscr F}_x,\,\Omega^k_x\otimes_{{\mathscr O}_x} {\mathscr F}_x \big)  \overset{\lambda_x}{\leftarrow}
	\Omega^k_x\otimes_{{\mathscr O}_x}{\rm End}_{{\mathscr O}_x}({\mathscr F}_x).$$
	As ${\mathscr F}$ is coherent, $\alpha_x$ is isomorphism \cite[p.\,239]{G-R}. Furthermore, if ${\mathscr F}_x $ is free, then $\lambda_x$ is bijective. 

If ${\mathcal S}$ is the singularity set of ${\mathscr F},$  it is also the singular locus of the sheaf 
 ${\mathscr Hom}_{{\mathscr O}_X}({\mathscr F},\,\Omega^2({\mathscr F})).$   
	Hence, for each point $x$ outside of ${\mathcal S},$  the fiber of ${\mathscr Hom}_{{\mathscr O}_X}({\mathscr F},\, \Omega^k({\mathscr F})) $ at  $x$
	can be identified with the vector space 
	$\Omega^k{(x)}\otimes{\rm End}({ F}{(x)}).$  
  According to (\ref{mathcalKnabla}),
	the curvature  ${\mathcal K}_{\nabla}$ of a holomorphic connection $\nabla$  determines the vector
	\begin{equation}\label{tildeKx}
	K_{\nabla}(x)\in \Omega^2(x)\otimes{\rm End}( F(x))
	\end{equation}
	for each $x\in X\setminus{\mathcal S}$. That is, $K_{\nabla}$ is a $2$-form ${\rm End}(F)$-valued. 
	
Let us assume that  ${\mathscr F}$ is  a Hermitian torsion-free sheaf on the K\"ahler manifold $X,$ 
 according to (\ref{(mathcalZ,Z)torsion}),
one defines 
\begin{equation}\label{Norm_Curvature}
\norm{ {\mathcal K}_{\nabla}}^2=({\mathcal K}_{\nabla},\,{\mathcal K}_{\nabla})=
 \int_{X\setminus{\mathcal N} } {\mathcal K}_{\nabla}\cdot {\mathcal K}_{\nabla} \,{\rm d}{\rm vol}=
\int_{X \setminus{\mathcal N}}   |K_{\nabla}\wedge\star\,   K_{\nabla}|,
\end{equation}
where  $|\,\cdot\,|$ is the corresponding norm on ${\rm End}(F)$ and $\star$ is the Hodge star operator.

 More concretely, if locally $  K_{\nabla}$ can be expressed as $\alpha \otimes A$, with $\alpha$ a $2$-form and $A$ a local section of ${\rm End}(F)$, then
the integrand in (\ref{Norm_Curvature}) is $(\alpha\wedge\star\alpha)\langle A\circ A\rangle$. In a local unitary frame of ${\rm End}(F)$, if the connection is compatilble  with the metric,  the matrix $\check A$ associated to $A$ is antihermitian and  
$\langle A\circ A\rangle=(-1/2) {\rm tr}(\check A\check A)$. That is,
\begin{equation}\label{traza-Norma}
 | K_{\nabla}\wedge\star\,   K_{\nabla}|=-(1/2){\rm tr}\big(  K_{\nabla}\wedge\star\,   K_{\nabla}\big).
\end{equation}

On the space of holomorphic gauge fields on the torsion-free sheaf ${\mathscr F}$ one can  define the following map 
$\,\mathcal{YM}:\nabla\in {\rm Hom}({\mathscr F},\,\Omega^1({\mathscr F}))\mapsto 
 \norm{ {\mathcal K}_{\nabla}}^2.\,$
 It is called the  Yang-Mills' functional. The $\nabla$ on which this functional 
takes a stationary value are the {\it holomorphic Yang-Mills fields}.

\subsubsection{Yang-Mills fields on reflexive sheaves}
If ${\mathscr G}$ is a reflexive sheaf  on $X,$  then it is torsion-free and  the codimension of the singular set  is  $\geq 3$
 \cite[Cap.\,V, Cor.\,5.20]{Koba},  \cite[Cor.\,1.4]{Hartshorne1}. Furthermore, if $C$ is a closed subset of $X$ with codimension  $\geq 2,$ then the restriction $\Gamma(X,\,{\mathscr G})\to \Gamma(X\setminus C,\,{\mathscr G})$ is an isomorphism \cite[Cap.\,V, Prop.\,5.21]{Koba} \cite[Prop.\,1.11]{Hartshorne2}.

    Let us suppose that ${\mathscr F}$ is a {\em reflexive} sheaf on the K\"ahler manifold $X$ endowed with a Hermitian metric. 
  Then ${\mathscr End}({\mathscr F})$ is also a reflexive sheaf \cite[Chap.\,V, Prop.\,4.15]{Koba}. If ${\mathcal S}$ is the singular locus of ${\mathscr F},$ 
   let  ${\mathcal W}$ be a closed subspace of $X$  containing ${\mathcal S}$ and such that ${\rm codim}\,{\mathcal W}=2.$ Let ${\mathcal N}$ denote a tubular neighborhood of ${\mathcal S}$ in ${\mathcal W}.$ Therefore,
  the vector space 
	$\Gamma(X\setminus {\mathcal N},\,\Omega^k  ({\mathscr End}({\mathscr F})))\simeq
	\Gamma(X,\,\Omega^k ({\mathscr End}({\mathscr F}))),$ and the latter is finite dimensional. 
	Furthermore, as  $(X\setminus{\mathcal N})\cap{\mathcal S}=\emptyset,\,$
 ${\mathscr Hom }_{{\mathscr O}_X}({\mathscr F},\, \Omega^k({\mathscr F}))|_{X\setminus{\mathcal N}}$ is a locally free sheaf. 
	 
	 Let $\nabla$ be a holomorphic  gauge field on ${\mathscr F}$, then 
	$${\mathcal K}_{\nabla} \in \Gamma\big(X,\,{\mathscr Hom }({\mathscr F},\, \Omega^2({\mathscr F}))\big)=\Gamma(X\setminus{\mathcal N},\, {\mathscr Hom }({\mathscr F},\, \Omega^2({\mathscr F}))\big)$$
	 can considered as a global section of  the   locally free sheaf
	 $\big(\Omega^2\otimes {\rm End}(F)\big)|_{X\setminus{\mathcal S}}.$
	 A similar observation is also valid for $\nabla.$

\smallskip

Let us assume that the  sheaf ${\mathscr F}$ supports a holomorphic connection   $\nabla_0$. By Proposition \ref{P:Atiyah}, given 
	${\mathcal E}_1,\dots {\mathcal E}_m$, a basis of 
	\begin{equation}\label{HOM-hom}
	{\rm Hom}({\mathscr F},\,\Omega^1({\mathscr F}))=
	\Gamma\big(Y,\, {\mathscr Hom}_{{\mathscr O}_X}({\mathscr F},\Omega^1({\mathscr F}))\big)= 
	\Gamma(Y\setminus{\mathcal S},\,\Omega^1\otimes{\rm End}( F )),
	\end{equation}
	 any holomorphic gauge field can be written
	$\nabla=\nabla_0+\sum\lambda_i{\mathcal E}_i,$ with $\lambda_i\in{\mathbb C}$. The curvature
	\begin{equation}\label{nabla-nabla0}
	{\mathcal K}_{\nabla}=\nabla\circ\nabla={\mathcal K}_{\nabla_0} +\sum_i\lambda_i {\mathcal B}_i+\sum_{ij}\lambda_i\lambda_j {\mathcal B}_{ij},
	\end{equation}
	where ${\mathcal B}_i:=\nabla_0({\mathcal E}_i)$ and ${\mathcal B}_{ij}:={\mathcal E}_i\wedge {\mathcal E}_j$.
	If moreover ${\mathscr F}$ hermitian, then   
	$$\norm{{\mathcal K}_{\nabla}  }^2=({\mathcal K}_{\nabla},\,{\mathcal K}_{\nabla})=P(\lambda_1,\dots,\lambda_m),$$
	where $P$ is a polynomial 
 of degree $\leq 4$ in the variables $\lambda_i$. 
	
	The Yang-Mills fields are those $\nabla$ defined by   constants $\lambda_i$ which satisfy the algebraic equations of degree $\leq 3$
 \begin{equation}\label{partialP}
 \frac{\partial\,P}{\partial \lambda_i}=0,\;\;\; i=1,\dots,m.
  \end{equation}
  Therefore,
	\begin{Thm}\label{P:numeroYM}
	 If the Hermitian reflexive sheaf ${\mathscr F}$  admits a holomorphic gauge field and 
	$m={\rm dim}\, {\rm Hom}({\mathscr F},\,\Omega^1({\mathscr F}) ).$
	Then the set ${\sf YM}({\mathscr F})$ of holomorphic Yang-Mills fields on ${\mathscr F}$ is in bijective correspondence with the points of the algebraic set in ${\mathbb C}^m$ defined by $m$  algebraic equations. 
	
	In particular, if 
	 $m=2$ and the cardinal of ${\sf YM}({\mathscr F})$ is finite, then
	$\# {\sf YM}({\mathscr F})  \leq 9$.  
	\end{Thm}
	{\it Proof.} The case $m=2$ is a consequence of B\'ezout's theorem. \qed

 \smallskip
 
 Taking into account Proposition \ref{P:Atiyah} and (\ref{HOM-hom}), any ``variation'' of a holomorphic connection $\nabla$ on ${\mathscr F}$ can be written as $\nabla_{\epsilon}=\nabla+\epsilon E,$
 with $\epsilon\in{\mathbb C}$ and $E\in\Gamma(X\setminus S,\,\Omega^1\otimes {\rm End}(F)).$
\begin{equation}\label{(1/2)}
(1/2)\frac{d}{d\epsilon}\Big|_{\epsilon=0}||K_{\nabla_{\epsilon}}||^2=\int_{X\setminus{\mathcal N}} \langle K_{\nabla},\, \nabla  E\rangle \, {\rm d vol}=: (  K_{\nabla},\, \nabla{ E}).
\end{equation}
Therefore, $\nabla$ is a Yang-Mills field if for any ``variation'' ${\mathcal E}$ of $\nabla$ 
\begin{equation}\label{YanMillsEq}
(  K_{\nabla},\, \nabla{ E})=0.
\end{equation}
 In particular, the flat holomorphic gauge fields are Yang-Mills.

	We denote 
	$${}^{(p)}\nabla:\Gamma(X,\,\Omega^p ({\mathscr End}({\mathscr F})))\to
	\Gamma(X,\,\Omega^{p+1}({\mathscr End}({\mathscr F}))),$$
	the operator defined by the connection $\nabla$. In this notation Bianchi's identity is read as 
	$$ {}^{(2)}\nabla K_{\nabla}=0.$$
	
	On the other hand,   the orthogonality condition (\ref{YanMillsEq}) which satisfy the Yang-Mills fields gives rise to the following proposition.
	\begin{Prop}\label{P:ortho}
	The   holomorphic gauge field $\nabla$ on the Hermitian reflexive sheaf ${\mathscr F}$ is a Yang-Mills field iff
	its curvature $K_{\nabla}\in \Gamma(Y,\,\Omega^2\otimes_{\mathscr O}{\rm End}(F))$ is orthogonal to the vector space ${\rm Im}\,(^{(1)}\nabla)$.
	\end{Prop}



	\subsubsection{The case   ${\rm rank}\,{\mathscr F}=1$}
	In this case  the reflexive sheaf ${\mathscr F}$ is a locally free sheaf \cite[Prop.\,1.9]{Hartshorne1}.
	 Then ${\mathscr End}({\mathscr F})$ is the sheaf associated to the trivial line bundle ${\mathbb C}\times X\to X$. If $s$ is a local frame of the corresponding line bundle $F,$ a 
	holomorphic connection $\nabla$ on ${\mathscr F}$  on this frame is determined by a ${\mathbb C}$-valued $1$-form $A$, $\nabla s=A s$. In this frame $^{(p)}\nabla(\beta)=\partial\beta+A\wedge\beta-(-1)^p\beta\wedge A=\partial\beta,$
	for any ${\mathbb C}$-valued $p$-form $\beta$. That is, 
	\begin{equation}\label{pnabla=partial}
	^{(p)}\nabla=\partial.
	\end{equation}
	
	 In this case (\ref{nabla-nabla0}) reduces to $K_{\nabla}=K_{\nabla_0}$, since $\nabla_0(E_i)=\partial E_i=0.$ In particular, the Yang-Mills functional is constant.
	 
	  On the other hand, the Bianchi's identity reduces 
	   to $\partial K_{\nabla}=0$. If furthermore, $\nabla$ is a Yang-Mills field, from (\ref{YanMillsEq}) it follows  $\partial^{\dagger}(K_{\nabla})=0$, where $\partial^{\dagger}$ is the adjoint of $\partial$.
	  That is, $K_{\nabla}$ is $\partial$-harmonic. As $X$ is a K\"ahler manifold,  $K_{\nabla}$ is also $d$-harmonic. Hence, the norm of $K_{\nabla}$ minimizes the corresponding norm  in its cohomology class. That is, denoting by ${\bf c}$ the cohomology class defined by $K_{\nabla},$
	\begin{equation}\label{minimum}
	\mathcal{YM}(\nabla)=\norm{K_{\nabla}}^2={\min}\{\norm{\beta}^2\,|\, \beta\in{\bf c}\}.
	\end{equation}	  
	  
	As we are assuming that ${\mathscr F}$ supports a holomorphic gauge field, the first Chern class $c_1({\mathscr F})=[(2\pi)^{-1}K_{\nabla}]$ vanishes \cite{Atiyah}. Consequently, ${\bf c}=0$ and $\mathcal{YM}$ is the functional zero.


 \subsection{Yang-Mills fields on $B$-branes}\label{Ss:YM-Brenes}
In this subsection we assume that  $X$ be a smooth projective variety. Then  any  object 
of $D^b(X)$ is isomorphic to a bounded complex consisting of locally free sheaves \cite[Sect.\,36.36]{Stacks}. Thus, we assume that the $B$-brane ${\mathscr F}^{\bullet}$ on $X$ is a bounded complex of locally free sheaves on $X.$ In this  case, the jet sheaf 
${\mathscr J}^1({\mathscr F}^{\bullet})$ is the complex of vector bundles
${\mathscr F}^{\bullet}\oplus \Omega^1({\mathscr F}^{\bullet})$ endowed with
 ${\mathscr O}_X$-module structure  (\ref{fsigmabullet}).

A morphism in $D^b(X)$ from ${\mathscr M}^{\bullet}$ to ${\mathscr N}^{\bullet}$ is an equivalence class of pairs $[(s,\,f)],$ where 
${\mathscr M}^{\bullet}\overset{s}{\leftarrow}   {\mathscr G}^{\bullet}\overset{f}{\rightarrow}{\mathscr N}^{\bullet},$ where $s$ is a quasi-isomorphism. One has the morphisms induced on the cohomologies
$${\mathscr H}^i({\mathscr M}^{\bullet})\overset{(s^i)^{-1}}{\rightarrow}
{\mathscr H}^i({\mathscr G}^{\bullet})\overset{f^i}{\rightarrow}  {\mathscr H}^i({\mathscr N}^{\bullet}).$$
  Other pair $(t,\,g)$ with
${\mathscr M}^{\bullet}\overset{t}{\leftarrow}   \tilde{\mathscr G}^{\bullet}\overset{g}{\rightarrow}{\mathscr N}^{\bullet}$ equivalent to $(s,\,f)$ \cite[p.\,149]{Ge-Ma} defines the same morphism between the cohomologies.

 Thus, a gauge field on ${\mathscr F}^{\bullet}$, i.e. a morphism $\psi\in{\rm Hom}_{D^b(X)}\big( {\mathscr F}^{\bullet},\,{\mathscr J}^1({\mathscr F}^{\bullet})   \big)$ satisfying $\pi\circ \psi={\rm id},$ determines a unique morphism  
$$\psi^i:{\mathscr H}^i({\mathscr F}^{\bullet})\to {\mathscr H}^i({\mathscr J}^1({\mathscr F}^{\bullet}))={\mathscr H}^i({\mathscr F}^{\bullet})\,\tilde\oplus \, \Omega^1({\mathscr H}^i({\mathscr F}^{\bullet})),$$
such that the composition with the projection 
${\mathscr H}^i({\mathscr J}^1({\mathscr F}^{\bullet}))\overset{\pi^i} \to {\mathscr H}^i({\mathscr F}^{\bullet})$  is the 
identity on ${\mathscr H}^i({\mathscr F}^{\bullet}).$ 
 That is, $\psi^i$ is a holomorphic connection on the sheaf ${\mathscr H}^i({\mathscr F}^{\bullet}).$
 

  We  set $ \eta^j$ for the morphism of abelian sheaves defined by the inclusion in the direct sum
 $\eta^j: {\mathscr H}^j({\mathscr F}^{\bullet})  \rightarrow    {\mathscr H}^j(  {\mathscr F}^{\bullet})\oplus  \Omega^1({\mathscr H}^j(  {\mathscr F}^{\bullet})).$
 Hence, $\pi^j(\psi^j-\eta^j)=0$, 
  and  thus $\psi^j-\eta^j$ defines a morphism of abelian sheaves 
  \begin{equation}\label{vartethaConnections}
	\vartheta^j:{\mathscr H}^j({\mathscr F}^{\bullet})\to{\Omega^1}( {\mathscr H}^j({\mathscr F}^{\bullet})),
	\end{equation}
which, by (\ref{fsigmabullet}), satisfies  the Leibniz rule. 
 That is, 
\begin{Prop}\label{P:Connvartheta} The gauge field $\psi$ on the brane ${\mathscr F}^{\bullet}$ determines on each ${\mathscr O}_X$-module 
 ${\mathscr H}^j({\mathscr F}^{\bullet})$  a holomorphic connection $\vartheta^j$.
\end{Prop}

\begin{Rems}\label{R:variation}
Let $\psi,\phi$ be two gauge fields on ${\mathscr F}$, we set 
$$\xi:=\phi-\psi\in{\rm Hom}_{D^b(X)}
\big({\mathscr F}^{\bullet},\,{\mathscr J}^1({\mathscr F}^{\bullet})\big).$$
 Thus, $\xi$ determines a well defined morphism of ${\mathscr O}_X$-modules between the cohomologies, 
$\xi^j:{\mathscr H}^j({\mathscr F}^{\bullet})\to{\mathscr H}^j({\mathscr J}^1({\mathscr F}^{\bullet})).$

We denote by $\vartheta^j$ and $\chi^j$ the  connections on ${\mathscr H}^j({\mathscr F}^{\bullet})$ determined by $ \psi$ and $ \phi$, respectively.
Since $ \xi^j=(\phi^j-\eta^j)-(\psi^j-\eta^j)$, it follows that $\xi^j=\vartheta^j-\chi^j$. In short, $\xi^j$ is the ``variation'' on the connection $\vartheta^j$ induced by the ``variation'' $\xi$ of the gauge field $\psi$. 
\end{Rems}

The result deduced in the following paragraph gives us a suggestion for the definition of the  Yang-Mills functional on the gauge fields on a brane over a K\"ahler manifold.


\subsubsection{An Euler-Poincar\'e mapping.}\label{Sss: ConnectionsCS}  
Let ${\mathscr A}$ be a coherent sheaf on the K\"ahler manifold $X,$ and $\alpha:{\mathscr A}\to \Omega^1({\mathscr A})$ a
 holomorphic connection on ${\mathscr A}$. Denoting by ${\mathcal S}_{\mathscr A}$ the singular set of ${\mathscr A}$,   on
 $Y\setminus{\mathcal S}_{\mathscr A}$ we define differential form 
$$\Phi({\mathscr A},\alpha):={\rm tr}\big(K_{\alpha}\wedge\star K_{\alpha} \big)\in 
\Gamma(X\setminus{\mathcal S}_{\mathscr A},\,\Omega^{\rm top}\big) ,$$ 
 where $K_{\alpha}$ is the  curvature of $\alpha$, considered as an  ${\rm End}(A)$-valued $2$-form. 

By $\mathfrak{C}$, we denote the category whose objects are pairs $({\mathscr A},\,\alpha)$. A morphism $f:({\mathscr A},\,\alpha)\to
({\mathscr B},\,\beta)$ is a morphism of coherent sheaves compatible with the connections; i.e. such that $(1\otimes f)\circ\alpha=\beta\circ  f$.
\begin{Prop}\label{P:EulerPoincsre}
If $0\to ({\mathscr A},\,\alpha)\overset{f}{\to} ({\mathscr B},\,\beta)\overset{g}{\to} ({\mathscr C},\,\gamma)\to 0,$
is an exact sequence in $\mathfrak{C}$, then on $X\setminus {\mathcal S}$
$$\Phi({\mathscr B},\beta)=  \Phi({\mathscr A},\alpha)+\Phi({\mathscr C},\gamma),$$ 
 where   ${\mathcal S}$ is the union of the singular sets of ${\mathscr A},$ ${\mathscr B},$ and ${\mathscr C}$.
\end{Prop}
{\it Proof.} Let $x_0\in X\setminus {\mathcal S}$. As the exact sequence splits locally on $X\setminus {\mathcal S}$, there exists an open neighborhood $U$ of $x_0$ such that $g|_U$, in the sequence of locally free modules   
 $0\to {\mathscr A}|_U\overset{f|_U}{\to} {\mathscr B}|_U\overset{g|_U}{\to} {\mathscr C}|_U\to 0,$
has a right inverse $h$.

Let $a$ be a frame for ${\mathscr A}|_U$, then $\alpha(a)={\sf A }\cdot a$, where ${\sf A}$ is a matrix of $1$-forms on $U$. Furthermore, $a$ can be chosen so that ${\sf A}(x_0)=0$. Similarly, let $c$ be a frame for ${\mathscr C}|_U$, then $\gamma(c)={\sf C}\cdot c$ and we choose $c$ so that ${\sf C}(x_0)=0$.  From the splitting, it follows that $\{f(a),\, h(c)\}$ is a frame for ${\mathscr B}|_U$.
By the compatibility of the connections with $f$ and $g$, 
$$\beta(f(a))=(1\otimes f)(\alpha(a))=(1\otimes f)({\sf A} \cdot a)= {\sf A}\cdot f(a).$$
On the other hand, $\beta(h(c))={\sf R} \cdot f(a)+{\sf S}\cdot h(c) $, with ${\sf R}$ and ${\sf S}$ matrices of $1$-forms. But,
$$ {\sf C}\cdot c=\gamma(c)=\gamma(gh(c))=(1\otimes g)(\beta(h(c))=(1\otimes g)\big({\sf R}\cdot f(a)+{\sf S} \cdot h(c)\big).$$
As $g\circ f=0$ and $g\circ h=1$, it follows that ${\sf C}={\sf S}.$ That is, the matrix of the connection $\beta$ in the frame $\{f(a),\, h(c)\}$  is \begin{equation}\label{matrix Kb}
{\sf M}:=\begin{pmatrix}
{\sf A} & {\sf R} \\
0 & {\sf C}
\end{pmatrix}
\end{equation}

Since ${\sf A}(x_0)=0$ and ${\sf C}(x_0)=0$, the matrix of $K_{\alpha}(x_0)$, of the curvature of $\alpha$ at the point $x_0$, is
 $d{\sf A}$. Analogously,  the matrix of $K_{\gamma}(x_0)$ is $d{\sf C}$. The one of $K_{\beta}(x_0)$ is the exterior derivative of 
(\ref{matrix Kb}), since ${\sf M}\wedge {\sf M}=0$ at $x_0$. Then 
\begin{align}
{\rm tr }\big( K_{\beta}(x_0) \wedge\star K_{\beta}(x_0)  \big)&={\rm tr}( d{\sf A}\wedge\star d{\sf A} ) +
{\rm tr}( d{\sf C}\wedge\star d{\sf C} ) \notag \\  
&={\rm tr }\big( K_{\alpha}(x_0) \wedge\star K_{\alpha}(x_0)  \big)+{\rm tr }\big( K_{\gamma}(x_0) \wedge\star K_{\gamma}(x_0)  \big).\notag
\end{align}
As $x_0$ is an arbitrary point of $X\setminus{\mathcal S}$, it follows the proposition.
\qed

\smallskip

Let   $({\mathscr G}^{\bullet},\,{\delta}^{\bullet})$ be a bounded complex of {coherent}  sheaves on the manifold $X$.
Let $\nabla^{\bullet}$ be a family of holomorphic connections, compatible with the operators $\delta^{\bullet}$.  
That is, $\nabla^i:{\mathscr G}^{i}\rightarrow \Omega^1({\mathscr G}^{i})$ is a holomorphic connection on the coherent sheaf ${\mathscr G}^i$ such that
 $(1\otimes \delta^i)\nabla^i=\nabla^{i+1} \delta^i.$
Hence, $\nabla^i({\mathscr Ker}({\delta}^i))\subset {\mathscr Ker}(1\otimes{\delta}^i)$ and a similar relation for the image
  ${\mathscr Im}(\delta^{i-1})$.
It follows that $\nabla^i$ induces a connection $\theta^i$ on the cohomology
$\theta^i:{\mathscr H}^i({\mathscr G}^{\bullet})\to \Omega^1({\mathscr H}^i({\mathscr G}^{\bullet})).$
 Obviously, the restrictions of $\nabla^i$  determine connections on ${\mathscr Ker}(\delta^i)$ and
 ${\mathscr Im}(\delta^{i+1})$, respectively. One has the exact sequence
 \begin{equation}\label{exact_seq:C}
 0\to \big({\mathscr Ker}(\delta^i),\,\nabla^i\big)\to\big({\mathscr G}^i,\,\nabla^i\big)\to \big({\mathscr Im}(\delta^{i}),\,\nabla^{i+1}\big)\to 0
 \end{equation}
 in the category $\mathfrak{C}$. Similarly, we have the exact sequence
 \begin{equation}\label{exact_seq:C1}
 0\to  \big({\mathscr Im}(\delta^{i-1}),\,\nabla^{i}\big)\to\big({\mathscr Ker}(\delta^i),\,\nabla^i\big)\to\big({\mathscr H}^i,\,\theta^i\big)\to  0.
 \end{equation}

\begin{Cor}\label{Coro:Euler}
Denoting with ${\mathcal S}$ the union of the singular sets of the sheaves ${\mathscr G}^i$, then 
$\,\sum_{i}(-1)^i{\rm tr}\big(K_{\nabla^i}\wedge\star K_{\nabla^i}\big)=
  \sum_{i}(-1)^i{\rm tr}\big(K_{\theta^i}\wedge\star K_{\theta^i}\big)$
  on $X\setminus {\mathcal S}.$
\end{Cor}
  {\it Proof.} From Proposition \ref{P:EulerPoincsre} together with (\ref{exact_seq:C}), it follows
		$$\Phi({\mathscr G}^i,\,\nabla^i)= \Phi({\mathscr Ker}(\delta^i),\,\nabla^i)+\Phi({\mathscr Im}(\delta^i),\,\nabla^{i+1}).$$
	From (\ref{exact_seq:C1}), one obtains an analogous relation. Taking the alternate sums 
  $$\sum_i(-1)^i\Phi({\mathscr G}^{i},\,\nabla^{i})=\sum_i(-1)^i\Phi({\mathscr H}^{i},\,\theta^{i}). \qed $$


	\subsubsection{The Yang-Mills functional}  
We propose a definition for the Yang-Mills functional over gauge fields on a brane. This proposal is based on the following considerations:
\begin{enumerate}
  \item It is reasonable to require that this definition generalizes the one for coherent sheaves.
  \item As a gauge field $\psi$ on ${\mathscr F}^{\bullet}$ is a class of ``roofs'' \cite[p.\,148]{Ge-Ma} from the complex 
  ${\mathscr F}^{\bullet}$ to ${\mathscr J}^1({\mathscr F}^{\bullet}),$  and equivalent roofs determine the same morphisms between the cohomologies,
   it seems convenient to move on  the cohomology of these complexes.
  \item Let $E^{\bullet}$ be a bounded complex of   Hermitian vector bundles over the K\"ahler manifold $Y,$ and $\nabla^{\bullet}$ a family of connections compatible with 
	the Hermitian metrics and 
	the coboundary operators. Denoting by ${\mathscr H}^i(E^{\bullet})$ the cohomology sheaves, there exist connections $\theta^i$ on those sheaves, induced by the family $\nabla^{\bullet}$. By Corollary \ref{Coro:Euler} together with   (\ref{traza-Norma}), one has
  the following equality of Euler-Poincar\'e type. 
  $$\sum_i(-1)^i\norm{K_{\nabla^i}}^2  =
  \sum_i(-1)^i\norm{(K_{\theta^i}}^2.$$ 
 \end{enumerate}
 On the basis of the above considerations, it seems appropriate to define the value of the Yang-Mills functional on the    gauge $\psi$ on the brane ${\mathscr F}^{\bullet}$ as
 $\sum_i(-1)^i\norm{K_{\vartheta_i}}^2.$
   More precisely, taking into account Proposition \ref{P:Connvartheta}, we adopt the following definitions.   

\begin{Def}\label{D:HermitianBrane}
The brane $({\mathscr F}^{\bullet},\,\delta^{\bullet})$ is called a Hermitian brane, if the cohomology sheaves ${\mathscr H}^j$ are Hermitian
  ${\mathscr O}_X$-modules.
\end{Def}

Let $({\mathscr F}^{\bullet},\,{\delta}^{\bullet})$ be a Hermitian brane on the K\"ahler manifold $X$. Given a gauge
 field $\psi$  
on the brane $({\mathscr F}^{\bullet},\,\delta^{\bullet})$,
by Proposition \ref{P:Connvartheta}, one has the family of curvatures ${\mathcal K}_{\vartheta^i}$ of the connections induced on the cohomologies,
whose norms $\norm{{\mathcal K}_{\vartheta^i}}$  are defined in accordance with (\ref{Norm_Curvature}).

\begin{Def}\label{Def:YM2}
 Given a gauge field $\psi$ on the Hermitian $B$-brane $({\mathscr F}^{\bullet},\,\delta^{\bullet})$, if the sheaves ${\mathscr H}^i({\mathscr F}^{\bullet})$
 are reflexive, we define the value of the Yang-Mills functional at $\psi$ by
 \begin{equation}\label{Y-M(Brana)}
 \mathcal{YM}(\psi)=\sum_i(-1)^i\norm{{\mathcal K}_{\vartheta_i}}^2.
 \end{equation}
\end{Def}
	Thus, $\mathcal{YM}(\psi)$ 
	is a kind of Euler characteristic  of the gauge field.
	The {\it Yang-Mills fields} on the brane 
	$({\mathscr F}^{\bullet},\,{\delta}^{\bullet})$ are 
	 the critical points of the functional $\psi\mapsto \mathcal{YM}(\psi).$
	
	Note that if $({\mathscr F}^{\bullet},\,\delta^{\bullet})$ is an acyclic complex, then the Yang-Mills functional for this complex is identically zero. 
	\smallskip
		
		\smallskip
	
Let ${\mathscr A}^{\bullet}:=\big( {\mathscr A}^{\bullet},\,d_A^{\bullet},\,\alpha^ {\bullet }\big)$ be a complex in the category ${\mathfrak C}$; i.e, a complex of coherent sheaves with a family of holomorphic connections compatible with the coboundary operator $d_A$. Let $f:=(f^{\bullet})$ a morphism 
  $f^{\bullet}:{\mathscr A}^{\bullet}\to{\mathscr B}^{\bullet}$ in   ${\mathfrak C}$; that is, $f$ is a morphism of complexes compatible with the connections.    Let ${\mathscr C}^{\bullet}$ denote the mapping cone of 
  $f$. Thus, 
${\mathscr C}^{\bullet}=\big( {\mathscr A}^{\bullet}[1] \oplus {\mathscr B}^{\bullet},d_C^{\bullet},\nabla^{\bullet},
 \big)$, with $d_C(a,\,b)=\big(d(a),\,(-1)^{{\rm degree}\,a}f(a)+db   \big)$ and $\nabla(a,\,b)=(\alpha(a),\,\beta(b)).$ In fact, $(1\otimes d_C)\circ\nabla=\nabla\circ d_C$ and thus ${\mathscr C}^{\bullet}$ is a complex of the category ${\mathfrak C}$
 
 For each $i$ one has the following exact sequence in the category ${\mathfrak C}$
 $$0\to{\mathscr B}^{i}\to {\mathscr C}^{i}\to {\mathscr A}^{i}[1]\to 0.$$
 From Proposition \ref{P:EulerPoincsre},
 $\Phi({\mathscr B}^{i})+ \Phi({\mathscr A}^{i+1})=\Phi({\mathscr C}^{i})$. Multiplying by $(-1)^i$ and summing
\begin{align} \label{alingsumi}\sum_i(-1)^i{\rm tr}(K_{\beta^i}\wedge\star K_{\beta^i}) + &
 \sum_i(-1)^i{\rm tr}(K_{\alpha^{i+1}}\wedge\star K_{\alpha^{i+1}}) \\ \notag
 &=\sum_i(-1)^i{\rm tr}(K_{\nabla^i}\wedge\star K_{\nabla^i}).
\end{align}

Let us assume that
\begin{itemize}
 \item ${\mathscr A}^{i}$ and ${\mathscr B}^{i}$   Hermitian sheaves for all $i$.
 \item  $\alpha^i$ and $\beta^i$  are Hermitian gauge fields (i.e., compatible with the metric) on ${\mathscr A}^{i}$ and ${\mathscr B}^{i}$, respectively.
 \end{itemize}
Then one defines on ${\mathscr C}^i$  the metric
$\langle(a,b),\,(a',b')\rangle:=\langle a,\,a'\rangle+\langle b,\,b'\rangle$. The connection $\nabla^i$ is compatible with this metric.
 From the  equality  (\ref{alingsumi})  together with (\ref{traza-Norma}), one deduces the following proposition. 
\begin{Prop}\label{P:YMCone}
 With the above notations and under the above hypotheses,  
  $\alpha$ and $\beta$ determine in a natural way a gauge field $\nabla$ on the mapping cone of $f^{\bullet}$ satisfying 
 \begin{equation}\label{YMCone}
  \mathcal{YM}(\beta)-\mathcal{YM}(\alpha) = \mathcal{YM}(\nabla).
  \end{equation}
\end{Prop}
On the other hand, in the context of the branes theory,
the fact that the branes ${\mathscr  A}^{\bullet},\,{\mathscr  B}^{\bullet}$ and ${\mathscr  C}^{\bullet}$  are the members of 
the distinguished triangle ${\mathscr A}^{\bullet}\to {\mathscr B}^{\bullet}\to 
{\mathscr C}^{\bullet}\overset{+1}{\to}$ means that 
${\mathscr A}^{\bullet}$ and  ${\mathscr C}^{\bullet}$   can potentially bind together to form the  membrane ${\mathscr B}^{\bullet}$ \cite[Section 6.2.1]{Aspin}. Thus, the additive nature of equation
(\ref{YMCone}) is consistent with this interpretation.

	
	\smallskip	
		
From now on, we assume that ${\mathscr F}^{\bullet}$ is a {\em Hermitian $B$-brane} such that the cohomology sheaves ${\mathscr H}^i({\mathscr F}^{\bullet})$ are {\em reflexive}.

Let us suppose that the set of gauge fields on the 
 brane ${\mathscr F}^{\bullet}$ is nonempty. Let $m$ be the dimension of the vector space ${\rm Ext}^0({\mathscr F}^{\bullet},\,\Omega^1({\mathscr F}^{\bullet})).$ We denote by $\xi_1,\dots, \xi_m$ a basis of this vector space. According to Proposition 
\ref{P:Atiyah}, any gauge field $\psi$ on the brane can be expressed 
$$\psi=\tilde\psi+\sum_a\lambda_a\xi_a$$
 $\tilde\psi$ being a fixed gauge field and $\lambda_a\in{\mathbb C}.$ 
 Hence, the connections on the cohomology sheaves ${\mathscr H}^i$ 
 can be written in the form (see Remark \ref{R:variation})
 $$\vartheta^i=\tilde\vartheta^i+\sum_ a\lambda_a \xi^i_a,$$
 with $\xi^i_a\in{\rm Hom}({\mathscr H}^i,\,\Omega^1({\mathscr H}^i))$.
 The corresponding curvatures satisfy
 $$K_{\vartheta^i}= K_{\tilde\vartheta^i}+\sum_a \lambda_a\tilde\vartheta^i(\xi^i_a)+\sum_{a,b}\lambda_a\lambda_b\xi^i_a\wedge\xi^i_b.$$
 Therefore, 
$\norm{K_{\vartheta^i}}^2$ is a   polynomial $P^i(\lambda_1,\dots,\lambda_m)$ of degree $\leq 4.$ Thus, the critical points of the Yang-Mills functional correspond to the points $(\lambda_1,\dots,\lambda_m)\in{\mathbb C}^m$ which satisfy the equations
$\frac{\partial P}{\partial\lambda _a}=0,$ where $P$ is the   polynomial 
 $\sum_i(-1)^iP^i(\lambda_1,\dots, \lambda_m)$. We have the following result, which generalizes Theorem \ref{P:numeroYM}. 
\begin{Thm}\label{P:numero YM}
Assumed the set of gauge fields on the brane ${\mathscr F}^{\bullet}$ is nonempty and  
$m={\rm dim}\,{\rm Ext}^0({\mathscr F}^{\bullet},\,\Omega^1({\mathscr F}^{\bullet}))$. Then the set of Yang-Mills fields on ${\mathscr F}^{\bullet}$ are in bijective correspondence with the points of a subvariety of ${\mathbb C}^m$ defined by $m$  polynomials of degree $\leq 3$.   
\end{Thm}

	\smallskip
Let  $({\mathscr F}^{\bullet},\,\delta^{\bullet})$ be a $B$-brane, such that the
${\mathscr H}^i({\mathscr F}^{\bullet})$ are reflexive sheaves.
	If $\phi$ and $\psi$ are
	gauge fields on ${\mathscr F}^{\bullet}$ and $\xi=\phi-\psi$, using the notations introduced in Remark \ref{R:variation},
	the connections on the cohomologies induced by $\phi$ and $\psi$ satisfy $\chi^j(\xi)=\vartheta^j+\xi^j$, with 
	$$\xi^j\in \Gamma(Y\setminus {\mathcal S},\,\Omega^1\otimes_{\mathscr O}{\mathscr End}(H^j)),$$
	${\mathcal S}$ being the union of the singularity sets of the reflexive sheaves ${\mathscr H}^j.$
	
	With the mentioned notation, 
	an infinitesimal  variation $\psi_{\epsilon}$ of $\psi$ is 
	an element of the form $\epsilon\xi,$ where $\epsilon\in{\mathbb C}$ and $\xi\in {\rm Ext}^0({\mathscr F}^{\bullet},\,\Omega^1({\mathscr F}^{\bullet}).$ The relation $\psi_{\epsilon}=\psi+\epsilon\xi$ gives rise to the following equality between the connections on the cohomologies
	$$\vartheta^j_{\epsilon}=\vartheta^j+\epsilon\xi^j.$$
	Furthermore, on $Y\setminus {\mathcal S}$ the curvatures satisfy
	$$K_{\vartheta^j_{\epsilon}}=K_{\vartheta^j}+\epsilon\vartheta^j({\xi^j})+O(\epsilon^2),$$
	$\vartheta^j({\xi^j})$ being the covariant derivative of ${\xi^j}$ considered as a section 
	of $\Omega^1\otimes_{\mathscr O}{\mathscr End}(H^j)$.
	
From (\ref{(1/2)}) together with Definition \ref{Def:YM2},	the functional $\mathcal{YM}$ takes at the  gauge field ${\psi}$ a stationary value if
	\begin{equation}\label{(1/2)bis}
	(1/2)\frac{d}{d\epsilon}\mathcal{YM}(\psi_{\epsilon} )\Big|_{\epsilon =0}=\sum_j(-1)^j\langle K_{\vartheta^j},\,\vartheta^j(\xi^j)\rangle=0,
	\end{equation}
	for all   any  variation of $\psi$.
	In particular, if $\vartheta^i$ es a Yang-Mills field for all $i,$  then by (\ref{YanMillsEq}) 
	$\langle K_{\vartheta^i},\,\vartheta^i(\xi^i)\rangle=0,$
  and we have the following proposition.
  
	\begin{Prop}\label{estationaryY-Mfunctional}
Let	  $\psi$  gauge field   on the brane $({\mathscr F}^{\bullet},\,\delta^{\bullet}).$
	If  $\vartheta^i$ is a Yang-Mills field
	on ${\mathscr H}^i$ for all $i$, then $\psi$
	is a Yang-Mills field on the brane.
	\end{Prop}

	The following theorem is a partial  converse  to  Proposition \ref{estationaryY-Mfunctional}. 
	\begin{Thm}\label{Th:inverse} 
	 If each ${\mathscr F}^i$ is a semisimple object in the category of coherent sheaves on $X$ 
	and $\nabla^{\bullet}$ is a Yang-Mills field on   ${\mathscr F}^{\bullet},$ then the connection $\vartheta^j$ induced on ${\mathscr H}^j$ is a Yang-Mills field on this sheaf. 
	\end{Thm}
 {\it Proof.} 
  Since ${\mathscr F}^i$ is semisimple the following short exact sequence of coherent ${\mathscr O}_X$-modules
 $$0 \to {\mathscr Ker}(\delta^i) \to {\mathscr F}^i\to {\mathscr Coim}(\delta^i)\to 0$$
 splits. That is, ${\mathscr F}^i\simeq {\mathscr Ker}(\delta^i) \oplus {\mathscr Coim}(\delta^i).$ Since ${\mathscr Ker}(\delta^i)$ is semisimple 
 the  exact sequence 
 $0\to {\mathscr Im}({\delta}^{i-1})\to {\mathscr Ker}(\delta^i) \to {\mathscr H}^i\to 0$ also splits
 Thus, 
 \begin{equation}\label{Fi=Hi+Gi}
 {\mathscr F}^i\simeq {\mathscr H}^i \oplus {\mathscr G}^i,
 \end{equation}
 where ${\mathscr G}^i$ is isomorphic to the direct sum of ${\mathscr Coim}(\delta^i)$ and ${\mathscr Im}(\delta^{i-1})$.
 
 On the other hand, the coboundary operator $\delta^i:{\mathscr F}^i\to {\mathscr F}^{i+1}$ induces via the isomorphisms  (\ref{Fi=Hi+Gi}) to the morphism
 \begin{equation}\label{Fi=Hi+Gi(1)}
  \delta^i:{\mathscr H}^i \oplus {\mathscr G}^i\to {\mathscr H}^{i+1} \oplus {\mathscr G}^{i+1},\;\;\;
  (a,\,b)\mapsto (0,\,\delta^ib).
  \end{equation}
 
 Given $\xi\in {\rm Hom}_{D^b(X)} \big({\mathscr F}^{\bullet},\,\Omega^1\otimes_{\mathscr O}{\mathscr F}^{\bullet}  \big)$, according to Remark   \ref{R:variation},  it determines $\xi^i\in {\rm Hom}\big({\mathscr H}^{i},\,\Omega^1\otimes_{\mathscr O}{\mathscr H}^{i}  \big).$  As $\nabla^{\bullet}$ is, by hypothesis, a Yang-Mills field then (\ref{(1/2)bis}) is satisfied. 
 
 Given $j$, a general ``variation'' of $\vartheta^j$ is defined by an element
 $\tau\in{\rm Hom}\big({\mathscr H}^j,\,\Omega^1\otimes_{\mathscr O}{\mathscr H}^j\big).$ Under the hypotheses of the proposition, we need to prove that
 $$\langle K_{\vartheta^j},\,\vartheta^j(\tau)\rangle=0,$$
 for any variation $\tau$.
 The morphism $\tau$ can be extended to a morphism
 $$C^i:{\mathscr H}^i\oplus{\mathscr G}^i\to \Omega^1\otimes_{\mathscr O}\big( {\mathscr H}^i\oplus{\mathscr G}^i  \big),$$
 where 
 $$C^i(a,\,b)=\begin{cases}(\tau(a),\,0),\;\;\text{if}\;\, i=j\\
 (0,\,0), \;\;\text{if}\;\, i\ne j
 \end{cases} $$
 Moreover, the $C^i$ are compatible with the coboundaries. For example for $i=j$, by  (\ref{Fi=Hi+Gi(1)}),
 $((1\otimes \delta^j)\circ { C}^j)(a,\,b)=(1\otimes\delta^j)(\tau(a),\,0)=0$; and 	${ C}^{j+1}\circ\delta^{j}(a,\,b)=0$. Thus, by the isomorphism (\ref{Fi=Hi+Gi}) the $C^i$ determine a morphism
 $\xi:{\mathscr F}^{\bullet}\to\Omega^1\otimes_{\mathscr O}{\mathscr F}^{\bullet}$ in the category $D^b(X)$, and the corresponding $\xi^i$ induced in the cohomologies are all $0$ except when $i=j$, in which case $\xi^j=\tau$. Hence,
  by (\ref{(1/2)bis})
	$$0= \sum_i(-1)^i\langle K_{\vartheta^i},\,\vartheta^i(\xi^i)\rangle  =(-1)^j \langle K_{{\vartheta}^j},\,\vartheta^j(\xi^j)\rangle=(-1)^j\langle K_{{\vartheta}^j},\,\vartheta^j(\tau)\rangle.$$
	This holds for any ``variation'' $\tau$ of $\vartheta^j$.
	That is, by (\ref{YanMillsEq}), $\vartheta^j$ is a Yang-Mills field on ${\mathscr H}^j$. \qed

\section{Final Remarks}
The following are some aspects and assumptions made in the article that could be considered to generalize or extend the  derived results.

\begin{Rems}
In the definition of the Yang-Mills functional on the gauge fields on brane ${\mathscr F}^{\bullet}$ we assumed that the cohomology sheaves ${\mathscr H}^j({\mathscr F}^{\bullet})$ are reflexive. It would be interesting to extend the definition to holomorphic gauge fields on an arbitrary $B$-brane. 
 \end{Rems}
 
 \begin{Rems}
 
 The semisimplicity assumption in Theorem \ref{Th:inverse}  is a strong hypothesis that would be desirable to weaken.
 
 \end{Rems}
 
 \begin{Rems}
 In Section \ref{S:Flag_varieties}, we considered the variety of complete flags in ${\mathbb C}^3$
  and proved 
  Theorem \ref{Thm:flag}, which provides an upper bound for the number of gauge fields on this variety. The restriction to this flag variety arises from the fact that, in this case, any Schubert subvariety 
$Z$ is smooth. Thus, from the vanishing of the first Betti number, we deduce that $H^{1,0}(Z)=0,$ 
  which is a necessary condition for Proposition \ref{Lr-LrO} to hold.
  Following an alternative approach, could it be possible to prove a version of Theorem \ref{Thm:flag} 
 for any complete flag variety?
  \end{Rems}

 \begin{Rems}

The definitions introduced in the article and the basic results pertain to $B$-branes over a complex manifold.
 One may wonder which parts of the article's content are generalizable to the context of branes over varieties or schemes. It is expected that anything admitting a formulation in categorical language, such as the definition of jet sheaf or Proposition \ref{P:Atiyah},  can be translated into the framework of scheme theory. However, the translation of other concepts originating in differential geometry, such as the Yang-Mills functional, is not straightforward.
 
 \end{Rems}


\end{document}